\documentclass[11pt,a4paper, fullpage]{amsart}

\usepackage[centertags]{amsmath}
\usepackage{latexsym,amsthm,amssymb,verbatim}
\usepackage{epsf}
\usepackage[scanall]{psfrag}
\usepackage{graphicx, color}

\newcommand{\appsection}[1]{\let\oldthesection\thesection
  \renewcommand{\thesection}{Appendix \oldthesection}
  \section{#1}\let\thesection\oldthesection}

\newtheorem{defin}{Definition}[section]
\newtheorem{prop-def}[defin]{Proposition-Definition}
\newtheorem{lem}[defin]{Lemma}
\newtheorem{thm}[defin]{Theorem}
\newtheorem{remark}[defin]{Remark}
\newtheorem{cor}[defin]{Corollary}

\newtheorem{conv}[defin]{Convention}
\newtheorem{ex}[defin]{Example}
\newtheorem{claim}{Claim}
\newcommand{\e}{ \hfill $\diamond$}

\begin{document}

\title[Algorithmic Problems]{Algorithmic Problems in Amalgams of Finite Groups}%
\author{L.Markus-Epstein}\footnote{Supported in part at the Technion by a fellowship of the Israel Council for Higher Education}%
\address{Department of Mathematics \\
Technion \\
Haifa 32000, Israel}%
\email{epstin@math.biu.ac.il}%
\maketitle

\begin{abstract}

Geometric methods proposed by Stallings \cite{stal} for treating
finitely generated subgroups of free groups  were successfully
used to solve a wide collection of decision problems for free
groups and their subgroups \cite{b-m-m-w, kap-m, mar_meak, m-s-w,
mvw, rvw, ventura}.

It turns out that Stallings' methods can be effectively
generalized for the class of amalgams of finite groups
\cite{m-foldings}. In the present paper we employ subgroup graphs
constructed by the generalized Stallings' folding algorithm,
presented in \cite{m-foldings}, to solve various algorithmic
problems in amalgams of finite groups.

\end{abstract}

\pagestyle{headings}
%


\section{Introduction} \label{section:Introduction}
Decision (or \emph{algorithmic}) problems is one of the classical
subjects of combinatorial group theory, originating in the three
fundamental decision problems  posed by Dehn \cite{dehn} in 1911:
the \emph{word problem} (which asks to answer whether a word over
the group generators represents the identity), the \emph{conjugacy
problem} (which asks to answer whether an arbitrary pair of words
over the group generators define conjugate elements) and the
\emph{isomorphism problem} (which asks to answer whether an
arbitrary pair of finite presentations determine isomorphic
groups).

Though  Dehn solved all three of these problems as restricted to
the canonical presentation of fundamental groups of closed
2-manifolds,  they are theoretically undecidable (unsolvable) in
general \cite{miller71, miller92}. However restrictions to some
particular classes of groups may yield surprisingly good results.
Remarkable examples include the solvability of the word problem in
one-relator groups (Magnus, see II.5.4 in \cite{l_s}) and in
hyperbolic groups (Gromov, see 2.3.B in \cite{gro}). The reader is
referred to the papers of  Miller \cite{miller71, miller92} for a
survey of decision problems for groups.

The groups considered in the present paper are   amalgams of
finite groups.  As is well known \cite{b-f}, these groups are
hyperbolic. Therefore the word problem in this class of groups is
solvable.
A natural generalization of the word problem is   the
\emph{(subgroup) membership problem} (or the \emph{generalized
word problem}), which asks to decide whether a word in the
generators of the group is an element of the given subgroup. An
efficient solution of the membership problem in amalgams of finite
groups was presented by the author in \cite{m-foldings}, where
graph theoretic methods for treating amalgams of finite groups
were developed. Namely, a finitely generated subgroup $H$ of an
amalgam $G=G_1 \ast_A G_2$ of finite groups is canonically
represented by a finite labelled graph $\Gamma(H)$. This graph
carries all the essential information about the subgroup $H$
itself, which enables one to ``read off'' a solution of the
membership problem in $H$ directly from its subgroup graph
$\Gamma(H)$. This yields a quadratic (and sometimes even linear)
time solution of the membership problem in amalgams of finite
groups.

Such strategy was originally developed by  Stallings \cite{stal}
to treat finitely generated subgroups of free groups. Stallings'
approach was topological. He showed that every finitely generated
subgroup of a free group is canonically represented by a minimal
immersion of a bouquet of circles.  Using the graph theoretic
language, the results of \cite{stal} can be restated as follows. A
finitely generated subgroup of a free group is canonically
represented by a finite labelled graph which can be constructed
algorithmically by a so called process of \emph{Stallings'
foldings} (\emph{Stallings' folding algorithm}). Moreover, this
algorithm is quadratic in the size of the input \cite{kap-m,
m-s-w}. See \cite{tuikan} for a faster implementation of this
algorithm.

In \cite{m-foldings}  Stallings' folding algorithm was generalized
to the class of amalgams of finite groups. Along the current paper
we refer to this algorithm as the \emph{generalized Stallings'
folding algorithm}. Its description is included in the Appendix.

Note that graphs constructed by the Stallings' folding algorithm
can be viewed as  finite inverse automata as well. This
convergence of ideas from the group theory, topology, graph
theory, the theory of finite automata and finite semigroups yields
reach computational and algorithmic results concerning free groups
and their subgroups. In particular, this approach gives polynomial
time algorithms to solve the membership problem, the finite index
problem, to compute closures of subgroups in various profinite
topologies. See \cite{b-m-m-w, mar_meak, m-s-w, mvw, rvw, ventura}
for these and other examples of the applications of the Stallings'
approach in free groups, and \cite{kap-w-m, kmrs, m_w, schupp} for
the applications in some other classes of groups. Note that the
Stallings' ideas were recast in a combinatorial graph theoretic
way in the remarkable survey paper of Kapovich and Myasnikov
\cite{kap-m}, where these methods were applied systematically to
study the subgroup structure of free groups.


Our  objective is to apply the generalized Stallings' methods
developed by the author in \cite{m-foldings} to solve various
decision problems concerning finitely generated subgroups of
amalgams of finite groups algorithmically (that is to find a
precise procedure, an \emph{algorithm}), which extends the results
of \cite{kap-m}.

Our results include polynomial solutions for the following
algorithmic problems in amalgams of finite groups, which are known
to be unsolvable in general \cite{miller71, miller92}:
\begin{itemize}
\item computing subgroup presentations,
%
\item detecting triviality of a given subgroup,
\item the freeness problem,
\item  the finite index problem,
\item the separability problem,
\item the conjugacy problem,
\item the normality,
\item the intersection problem,
\item the malnormality problem,
\item the power problem,
\item reading off Kurosh decomposition for finitely generated
subgroups of free products of finite groups.
\end{itemize}

These results are spread out between three papers: \cite{m-algII,
m-kurosh} and the current one. In \cite{m-kurosh} free products of
finite groups are considered, and an efficient procedure to read
off a Kurosh decomposition is presented.

The splitting between \cite{m-algII} and the current paper was
done with the following idea in mind. It turn out that some
subgroup properties, such as computing of a subgroup presentation
and index, as well as detecting of freeness and normality, can be
obtained directly by an analysis of the corresponding subgroup
graph.
 Solutions of others require some additional
constructions. Thus, for example, intersection properties can be
examined via product graphs, and separability needs constructions
of a pushout of graphs.

In the current paper  algorithmic problems of the first type are
presented: the computing of subgroup presentations, the freeness
problem and the finite index problem. The separability problem is
also included here, because it is closely related with the other
problems presented in the current paper. The rest of the
algorithmic problems are introduced in  \cite{m-algII}.


The paper is organized as follows. The Preliminary Section
includes the description of the basic notions  used along the
present paper. Readers familiar with  amalgams, normal words in
amalgams and labelled graphs can skip it. The next section
presents a summary of the results from \cite{m-foldings} which are
essential for our algorithmic purposes. It describes the nature
and the properties of the subgroup graphs constructed by the
generalized Stallings' folding  algorithm in \cite{m-foldings}.
The rest of the sections are titled by the names of various
algorithmic problems and present definitions (descriptions) and
solutions of the corresponding algorithmic problems. The relevant
references to other papers considering similar problems and a
rough analysis of the complexity of the presented solutions
(algorithms) are provided. In contrast with papers that establish
the exploration of the complexity of decision problems as their
main goal (for instance, \cite{generic-case, average-case,
tuikan}), we do it rapidly (sketchy) viewing in its analysis a way
to emphasize the effectiveness  of our methods.

\subsection*{Other Methods}  \
There have been a number of papers, where methods, not based on
Stallings' foldings, have been presented. One can use these
methods to treat finitely generated subgroups of amalgams of
finite groups. A topological approach can be found in works of
Bogopolskii \cite{b1, b2}. For the automata theoretic approach,
see papers of Holt and Hurt \cite{holt-decision, holt-hurt},
papers of Cremanns, Kuhn, Madlener and Otto \cite{c-otto,
k-m-otto}, as well as the recent paper of Lohrey and Senizergues
\cite{l-s}.

However the methods for treating finitely generated subgroups
presented in the above papers were applied to some particular
subgroup property. No one of these papers has as its goal a
solution of various algorithmic problems, which we consider as our
primary aim.  Moreover, similarly to the case of free groups (see
\cite{kap-m}), our combinatorial approach seems to be the most
natural one for this purpose.


\section{Acknowledgments}

I wish to deeply thank to my PhD advisor Prof. Stuart W. Margolis
for introducing me to this subject, for his help  and
encouragement throughout  my work on the thesis. I owe gratitude
to Prof. Arye Juhasz for his suggestions and many useful comments
during the writing of this paper. I gratefully acknowledge a
partial support at the Technion by a fellowship of the Israel
Council for Higher Education.

%
\section{Preliminaries} \label{section:Preliminaries}
\subsection*{Amalgams}

Let $G=G_1 \ast_{A} G_2$ be a free product of $G_1$ and $G_2$ with
amalgamation, customary, an \emph{amalgam} of $G_1$ and $G_2$.
We assume that the (free) factors are given by the  finite group
presentations
\begin{align} G_1=gp\langle X_1|R_1\rangle, \ \ G_2=gp\langle
X_2|R_2\rangle \ \ {\rm such \ that} \ \ X_1^{\pm} \cap
X_2^{\pm}=\emptyset. \tag{\text{$1.a$}}
\end{align}
 $A= \langle Y   \rangle$ is a group such that there exist two
monomorphisms
\begin{align}
\phi_1:A \rightarrow G_1 \ {\rm and } \ \phi_2:A \rightarrow G_2.
\tag{\text{$1.b$}}
\end{align}
Thus $G$ has a finite group presentation
\begin{align}
G=gp\langle X_1,X_2 | R_1, R_2, \phi_1(a)=\phi_2(a), \; a \in Y
\rangle. \tag{\text{$1.c$}}
\end{align}

We  put $X=X_1 \cup X_2$,  $R=R_1 \cup R_2 \cup
\{\phi_1(a)=\phi_2(a) \; | \; a \in Y \} $. Thus $G=gp\langle
X|R\rangle$.

As is well known \cite{l_s, m-k-s, serre}, the free factors embed
in $G$. It enables us to identify $A$ with its monomorphic image
in each one of the free factors. Sometimes in order to make the
context clear we use \fbox{$G_i \cap A$}
\footnote{Boxes are used for emphasizing purposes only.}
to denote the monomorphic image of $A$ in $G_i$ ($i \in \{1,2\}$).

Elements of $G=gp \langle X |R \rangle$ are equivalence classes of
words. However it is customary to blur the distinction between a
word $u$ and the equivalence class containing $u$. We will
distinguish between  them by using different equality signs:
\fbox{``$\equiv$''} for the equality of two words and
\fbox{``$=_G$''} to denote the equality of two elements of $G$,
that is the equality of two equivalence classes. Thus in
$G=gp\langle x \; | \; x^4 \rangle$, for example, $x \equiv x$ but
$x \not\equiv x^{-3}$, while $x=_G x^{-3}$.



\subsection*{Normal Forms}
Let $G=G_1 \ast_A G_2$. A word $g_1g_2 \cdots g_n \in G$ is
\emph{in normal form} (or, simply, it is a \emph{normal word}) if:
\begin{enumerate}
    \item [(1)] $g_i \neq_G 1$ lies in one of the  factors, $G_1$ or $G_2$,
    \item [(2)] $g_i$ and $g_{i+1}$ are in different factors,
    \item [(3)] if $n \neq 1$, then $g_i \not\in A$.
\end{enumerate}
We call the sequence $(g_1, g_2, \ldots, g_n)$ a
 \emph{normal decomposition} of the element $g \in  G $, where $g=_G g_1g_2 \cdots g_n$.

Any $g \in G$ has a representative in a normal form (see, for
instance, p.187 in \cite{l_s}).  If $g \equiv g_1g_2 \cdots g_n $
is in normal form and $n>1$, then the Normal Form Theorem (IV.2.6
in \cite{l_s}) implies that $g \neq_G 1$. The number $n$ is unique
for a given element $g$ of $G$ and it is called the \emph{syllable
length} of $g$. We denote it $l(g)$. We use $|g|$ to denote the
length of $g$ as a word in $X^*$.


\subsection*{Labelled graphs}
Below we follow the notation of \cite{gi_sep, stal}.

A graph $\Gamma$ consists of two sets $E(\Gamma)$ and $V(\Gamma)$,
and two functions $E(\Gamma)\rightarrow E(\Gamma)$  and
$E(\Gamma)\rightarrow V(\Gamma)$: for each $e \in E$ there is an
element $\overline{e} \in E(\Gamma)$ and an element $\iota(e) \in
V(\Gamma)$, such that $\overline{\overline{e}}=e$ and
$\overline{e} \neq e$.

The elements of $E(\Gamma)$ are called \textit{edges}, and an $e
\in E(\Gamma)$ is a \emph{direct edge} of $\Gamma$, $\overline{e}$
is the \emph{reverse (inverse) edge} of $e$.

The elements of $V(\Gamma)$ are called \textit{vertices},
$\iota(e)$ is the \emph{initial vertex} of $e$, and
$\tau(e)=\iota(\overline{e})$ is the \emph{terminal vertex} of
$e$. We call them the \emph{endpoints} of the edge $e$.

A  \emph{path of length $n$} is  a sequence of $n$ edges $p=e_1
\cdots  e_n $ such that $v_i=\tau(e_i)=\iota(e_{i+1})$ ($1 \leq
i<n$).  We call $p$ a \emph{path from $v_0=\iota(e_1)$ to
$v_n=\tau(e_n)$}. The \emph{inverse} of the path $p$ is
$\overline{p}=\overline{e_n} \cdots \overline{e_1}$. A path of
length 0 is the \emph{empty path}.

We say that the graph $\Gamma$ is \emph{connected} if $V(\Gamma)
\neq \emptyset$ and any two vertices  are joined by a path. The
path $p$ is \emph{closed} if $\iota(p)=\tau(p)$, and it is
\emph{freely reduced} if $e_{i+1} \neq \overline{e_i}$ ($1 \leq i
<n$). $\Gamma$ is a \emph{tree} if it is a connected graph and
every closed freely reduced path in $\Gamma$ is empty.

A \emph{subgraph} of $\Gamma$ is a graph $C$  such that $V(C)
\subseteq V(\Gamma)$ and $E(C) \subseteq E(\Gamma)$. In this case,
by abuse of language, we write $C\subseteq \Gamma$.
Similarly, whenever we write $\Gamma_1 \cup \Gamma_2$ or $\Gamma_1
\cap \Gamma_2$,  we always mean that the set operations are, in
fact,  applied to the vertex sets and the edge sets of the
corresponding graphs.


A \emph{labelling} of $\Gamma$ by the set $X^{\pm}$ is a function
$$lab: \: E(\Gamma)\rightarrow X^{\pm}$$ such that for each $e \in
E(\Gamma)$, $lab(\overline{e}) \equiv (lab(e))^{-1}$.

The last equality enables one, when representing the labelled
graph $\Gamma$ as a directed diagram,  to represent only
$X$-labelled edges, because $X^{-1}$-labelled edges can be deduced
immediately from them.

A graph with a labelling function is called a \emph{labelled (with
$X^{\pm}$) graph}.  The only graphs considered in the present
paper are labelled graphs.

A labelled graph is called \emph{well-labelled} if
$$\iota(e_1)=\iota(e_2), \; lab(e_1) \equiv lab(e_2)\ \Rightarrow \
e_1=e_2,$$ for each pair of edges $e_1, e_2 \in E(\Gamma)$. See
Figure \ref{fig: labelled, well-labelled graphs}.

\begin{figure}[!h]
\psfrag{a }[][]{$a$} \psfrag{b }[][]{$b$} \psfrag{c }[][]{$c$}
\psfrag{e }[][]{$e_1$}
\psfrag{f }[][]{$e_2$}
\psfragscanon \psfrag{G }[][]{{\Large $\Gamma_1$}}
\psfragscanon \psfrag{H }[][]{{\Large $\Gamma_2$}}
\psfragscanon \psfrag{K }[][]{{\Large $\Gamma_3$}}
\includegraphics[width=\textwidth]{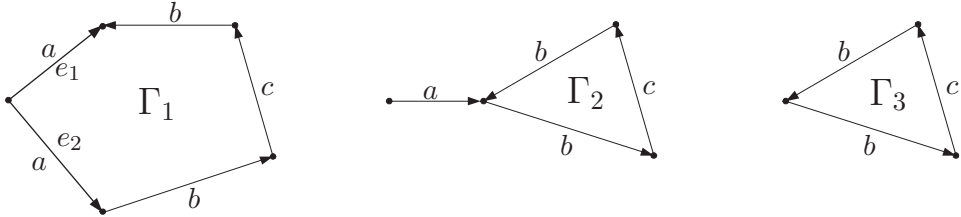}
\caption[The construction of $\Gamma(H_1)$]{ \footnotesize {The
graph $\Gamma_1$ is labelled with $\{a,b,c\} ^{\pm}$, but it is
not well-labelled. The graphs $\Gamma_2$ and $\Gamma_3$ are
well-labelled with $\{a,b,c\} ^{\pm}$.}
 \label{fig: labelled, well-labelled graphs}}
\end{figure}

If a finite graph $\Gamma$ is not well-labelled then a process of
iterative identifications of each pair  $\{e_1,e_2\}$ of distinct
edges with the same initial vertex and the same label to a single
edge yields a well-labelled graph. Such identifications are called
\emph{foldings}, and the whole process is known as the process of
\emph{Stallings' foldings} \cite{b-m-m-w, kap-m, mar_meak, m-s-w}.

 Thus the graph $\Gamma_2$ on Figure
\ref{fig: labelled, well-labelled graphs}  is obtained from the
graph $\Gamma_1$ by folding the edges $e_1$ and $e_2$ to a single
edge labelled by $a$.

Notice that the graph $\Gamma_3$ is obtained from the graph
$\Gamma_2$ by removing the edge labelled by $a$ whose initial
vertex has degree 1. Such an edge is called a \emph{hair}, and the
above procedure is used to be called \emph{``cutting hairs''}.


The label of a path $p=e_1e_2 \cdots e_n$ in $\Gamma$, where $e_i
\in E(\Gamma)$, is the word $$lab(p) \equiv lab(e_1)\cdots
lab(e_n) \in (X^{\pm})^*.$$ Notice that the label of the empty
path is the empty word. As usual, we identify the word $lab(p)$
with the corresponding element in $G=gp\langle X | R \rangle$. We
say that $p$ is   a \emph{normal path} (or $p$ is a path in
\emph{normal form}) if $lab(p)$ is a normal word.

If $\Gamma$ is a well-labelled graph then a path $p$ in $\Gamma$
is freely reduced if and only if $lab(p)$ is a freely reduced
word.
Otherwise $p$  can be converted into a freely reduced path $p'$ by
iterative  removals of the subpaths   $e\overline{e}$
(\emph{backtrackings}) (\cite{mar_meak, kap-m}).  Thus
$$\iota(p')=\iota(p), \ \tau(p')=\tau(p) \ \; {\rm and } \ \; lab(p)=_{FG(X)} lab(p'),$$ where \fbox{$FG(X)$} is a free group
with a free basis $X$. We say that $p'$ is obtained from $p$ by
\emph{free reductions}.


If $v_1,v_2 \in V(\Gamma)$ and $p$ is a path in $\Gamma$ such that
$$\iota(p)=v_1, \ \tau(p)=v_2 \ {\rm and } \ lab(p)\equiv u,$$
then, following the automata theoretic notation, we simply write
\fbox{$v_1 \cdot u=v_2$} to summarize this situation, and  say
that the word $u$ is \emph{readable} at $v_1$ in $\Gamma$.

A pair \fbox{$(\Gamma, v_0)$}  consisting  of the graph $\Gamma$
and the \emph{basepoint} $v_0$ (a distinguished vertex of the
graph $\Gamma$) is called  a \emph{pointed graph}.

Following the notation of Gitik (\cite{gi_sep}) we denote the set
of all closed paths in $\Gamma$ starting at $v_0$  by
\fbox{$Loop(\Gamma, v_0)$},  and the image of $lab(Loop(\Gamma,
v_0))$ in $G=gp\langle X | R \rangle$  by \fbox{$Lab(\Gamma,
v_0)$}. More precisely,
$$Loop(\Gamma, v_0)=\{ p \; | \; p {\rm  \ is \ a \ path \ in \ \Gamma \
with} \ \iota(p)=\tau(p)=v_0\}, $$
$$Lab(\Gamma,v_0)=\{g \in G \; | \; \exists p \in Loop(\Gamma,
v_0) \; : \; lab(p)=_G g \}.$$


It is easy to see that $Lab(\Gamma, v_0)$ is a subgroup of $G$
(\cite{gi_sep}). Moreover, $Lab(\Gamma,v)=gLab(\Gamma,u)g^{-1}$,
where $g=_G lab(p)$, and $p$ is a path in $\Gamma$ from $v$ to $u$
(\cite{kap-m}).
%
%
%
If $V(\Gamma)=\{v_0\}$ and $E(\Gamma)=\emptyset$ then we assume
that $H=\{1\}$.



We say that $H=Lab(\Gamma, v_0)$ is \emph{the subgroup of $G$
determined by the graph $(\Gamma,v_0)$}. Thus any pointed graph
labelled by $X^{\pm}$, where $X$ is a generating set of a group
$G$, determines a subgroup of $G$. This argues the use of the name
\emph{subgroup graphs} for such graphs.

\subsection*{Morphisms of Labelled Graphs} \label{sec:Morphisms
Of Well-Labelled Graphs}

Let $\Gamma$ and $\Delta$ be graphs labelled with $X^{\pm}$. The
map $\pi:\Gamma \rightarrow \Delta$ is called a \emph{morphism of
labelled graphs}, if $\pi$ takes vertices to vertices, edges to
edges, preserves labels of direct edges and has the property that
$$ \iota(\pi(e))=\pi(\iota(e)) \ {\rm and } \
\tau(\pi(e))=\pi(\tau(e)), \ \forall e\in E(\Gamma).$$
An injective morphism of labelled graphs is called an
\emph{embedding}. If $\pi$ is an embedding then we say that the
graph $\Gamma$ \emph{embeds} in the graph $\Delta$.


A \emph{morphism of pointed labelled graphs} $\pi:(\Gamma_1,v_1)
\rightarrow (\Gamma_2,v_2)$  is a morphism of underlying labelled
graphs $ \pi: \Gamma_1\rightarrow \Gamma_2$ which preserves the
basepoint $\pi(v_1)=v_2$. If $\Gamma_2$ is well-labelled then
there exists at most one such morphism (\cite{kap-m}).


\begin{remark}[\cite{kap-m}] \label{unique isomorphism}
{\rm  If two pointed well-labelled (with $X^{\pm}$) graphs
$(\Gamma_1,v_1)$ and $(\Gamma_2,v_2)$  are isomorphic, then there
exists a unique isomorphism $\pi:(\Gamma_1,v_1) \rightarrow
(\Gamma_2,v_2)$. Therefore $(\Gamma_1,v_1)$ and $(\Gamma_2,v_2)$
can be identified via $\pi$. In this case we sometimes write
$(\Gamma_1,v_1)=(\Gamma_2,v_2)$.} \e
\end{remark}

The notation $\Gamma_1=\Gamma_2$ means that there exists an
isomorphism between these two graphs. More precisely, one can find
$v_i \in V(\Gamma_i)$ ($i \in \{1,2\}$) such that
$(\Gamma_1,v_1)=(\Gamma_2,v_2)$ in the sense of Remark~\ref{unique
isomorphism}.


\section{Subgroup Graphs}

The current section is devoted to the discussion on subgroup
graphs constructed by the generalized Stallings' folding
algorithm. The main results of \cite{m-foldings} concerning these
graphs (more precisely, Theorem 7.1, Lemma 8.6, Lemma 8.7, Theorem
8.9 and Corollary 8.11 in \cite{m-foldings}), which are essential
for the current paper, are summarized in Theorem~\ref{thm:
properties of subgroup graphs} below. All the missing notations
are explained along the rest of the present section.



\begin{thm} \label{thm: properties of subgroup graphs}
Let $H=\langle h_1, \cdots, h_k \rangle$ be a finitely generated
subgroup of an amalgam of finite groups $ G=G_1 \ast_A G_2$.

Then there is an algorithm (\underline{the generalized Stallings'
folding  algorithm}) which  constructs a finite labelled graph
$(\Gamma(H),v_0)$ with the following properties:
\begin{itemize}
\item[(1)] $ {Lab(\Gamma(H),v_0)}= {H}. $

\item[(2)] Up to isomorphism, $(\Gamma(H),v_0)$ is  a unique
\underline{reduced precover} of $G$ determining $H$.

\item[(3)] A {\underline{normal word}} $g \in G$ is in $H$ if and
only if it labels a closed path in $\Gamma(H)$ starting at $v_0$,
that is $v_0 \cdot g=v_0$.

\item[(4)] Let $m$  be the sum of the lengths of words $h_1,
\ldots h_n$. Then the algorithm computes $(\Gamma(H),v_0)$ in time
$O(m^2)$.
Moreover, $|V(\Gamma(H))|$ and  $|E(\Gamma(H))|$ are proportional
to $m$.

\end{itemize}
\end{thm}

\begin{cor}
Theorem~\ref{thm: properties of subgroup graphs} (3) provides a
solution of the \underline{membership problem} for finitely
generated subgroups of amalgams of finite groups.
\end{cor}

Throughout the present paper the notation \fbox{$(\Gamma(H),v_0)$}
is  used for the finite labelled graph  constructed by the
generalized Stallings' folding  algorithm for a finitely generated
subgroup $H$ of an amalgam of finite groups $G=G_1 \ast_A G_2$.
%



\subsection*{Definition of Precovers:} The notion of
\emph{precovers} was defined by Gitik in \cite{gi_sep} for
subgroup graphs of amalgams. Below we present its definition and
list some  basic properties. In doing so, we rely on the notation
and results obtained in \cite{gi_sep}.

In \cite{m-foldings} some special cases of precovers,
\emph{reduced precovers}, were considered. However the properties
of \emph{reduced precovers} are irrelevant for the results
presented in the current paper. Hence we skip the discussion on
them, which can be found in \cite{m-foldings}.

Let $\Gamma$ be a graph labelled with $X^{\pm}$, where $X=X_1 \cup
X_2$ is the generating set of $G=G_1 \ast_A G_2$  given by
(1.a)-(1.c).
We view $\Gamma$ as a two colored graph: one color for each one of
the generating sets $X_1$ and $X_2$ of the factors $G_1$ and
$G_2$, respectively.

The vertex $v \in V(\Gamma)$ is called \emph{$X_i$-monochromatic}
if all the edges of $\Gamma$ incident with $v$ are labelled with
$X_i^{\pm}$, for some $i \in \{1,2\}$. We denote the set of
$X_i$-monochromatic vertices of $\Gamma$ by $VM_i(\Gamma)$ and put
$VM(\Gamma)= VM_1(\Gamma) \cup VM_2(\Gamma)$.

We say that a vertex $v \in V(\Gamma)$ is \emph{bichromatic} if
there exist edges $e_1$ and $e_2$ in $\Gamma$ with
$$\iota(e_1)=\iota(e_2)=v \ {\rm and} \ lab(e_i) \in X_i^{\pm}, \ i \in \{1,2\}.$$
The  set of bichromatic vertices of $\Gamma$ is denoted by
$VB(\Gamma)$.

A subgraph of $\Gamma$ is called \emph{monochromatic} if it is
labelled only with $X_1^{\pm}$ or only with $X_2^{\pm}$. An
\emph{$X_i$-monochromatic component} of $\Gamma$ ($i \in \{1,2\}$)
is a maximal connected subgraph of $\Gamma$ labelled with
$X_i^{\pm}$, which contains at least one edge.
Thus monochromatic components of $\Gamma$ are graphs determining
subgroups of the factors, $G_1$ or $G_2$.

We say that a graph $\Gamma$ is \emph{$ G$-based} if any path $p
\subseteq \Gamma$ with $lab(p)=_G 1$ is closed. Thus if $\Gamma$
is $G$-based then, obviously, it is well-labelled with $X^{\pm}$.

\begin{defin}[Definition of Precover] A $G$-based  graph $\Gamma$
is a \emph{precover} of $G$ if each $X_i$-monochromatic
component of $\Gamma$ is a \emph{cover} of $G_i$  ($i \in
\{1,2\}$).
\end{defin}

Following the terminology of Gitik (\cite{gi_sep}), we use the
term \emph{``covers of $G$''} for \emph{relative (coset) Cayley
graphs} of $G$ and denote by \fbox{$Cayley(G,S)$} the coset Cayley
graph of $G$ relative to the subgroup $S$ of
$G$.\footnote{Whenever the notation $Cayley(G,S)$ is used, it
always means that $S$ is a subgroup of the group $G$ and the
presentation of $G$ is fixed and clear from the context. }
If $S=\{1\}$, then $Cayley(G,S)$ is the \emph{Cayley graph} of $G$
and the notation \fbox{$Cayley(G)$} is used.

Note that the use of the term ``covers'' is adjusted by the  well
known fact that a geometric realization of a coset Cayley graph of
$G$ relative to some $S \leq G$ is a 1-skeleton of a topological
cover corresponding to $S$ of the standard 2-complex representing
the group $G$ (see \cite{stil}, pp.162-163).

\begin{conv}
By the above definition, a precover doesn't have to be a connected
graph. However along this paper we restrict our attention only to
connected precovers. Thus any time this term
 is used, we always mean that the corresponding graph
is connected unless it is stated otherwise.

We follow the convention that a graph $\Gamma$ with
$V(\Gamma)=\{v\}$ and $E(\Gamma)=\emptyset$ determining the
trivial subgroup (that is $Lab(\Gamma,v)=\{1\}$) is a (an empty)
precover of $G$.  \e
\end{conv}

\begin{figure}[!h]
\psfrag{x }[][]{$x$} \psfrag{y }[][]{$y$} \psfrag{v }[][]{$v$}
\psfrag{u }[][]{$u$}
\psfrag{w }[][]{$w$}
\psfrag{x1 - monochromatic vertex }[][]{{\footnotesize
$\{x\}$-monochromatic vertex}}
\psfrag{y1 - monochromatic vertex }[][]{\footnotesize
{$\{y\}$-monochromatic vertex}}
\psfrag{ bichromatic vertex }[][]{\footnotesize {bichromatic
vertex}}
\psfragscanon \psfrag{G }[][]{{\Large $\Gamma_1$}}
\psfragscanon \psfrag{K }[][]{{\Large $\Gamma_2$}}
\psfragscanon \psfrag{H }[][]{{\Large $\Gamma_3$}}
\psfragscanon \psfrag{L }[][]{{\Large $\Gamma_4$}}
\includegraphics[width=\textwidth]{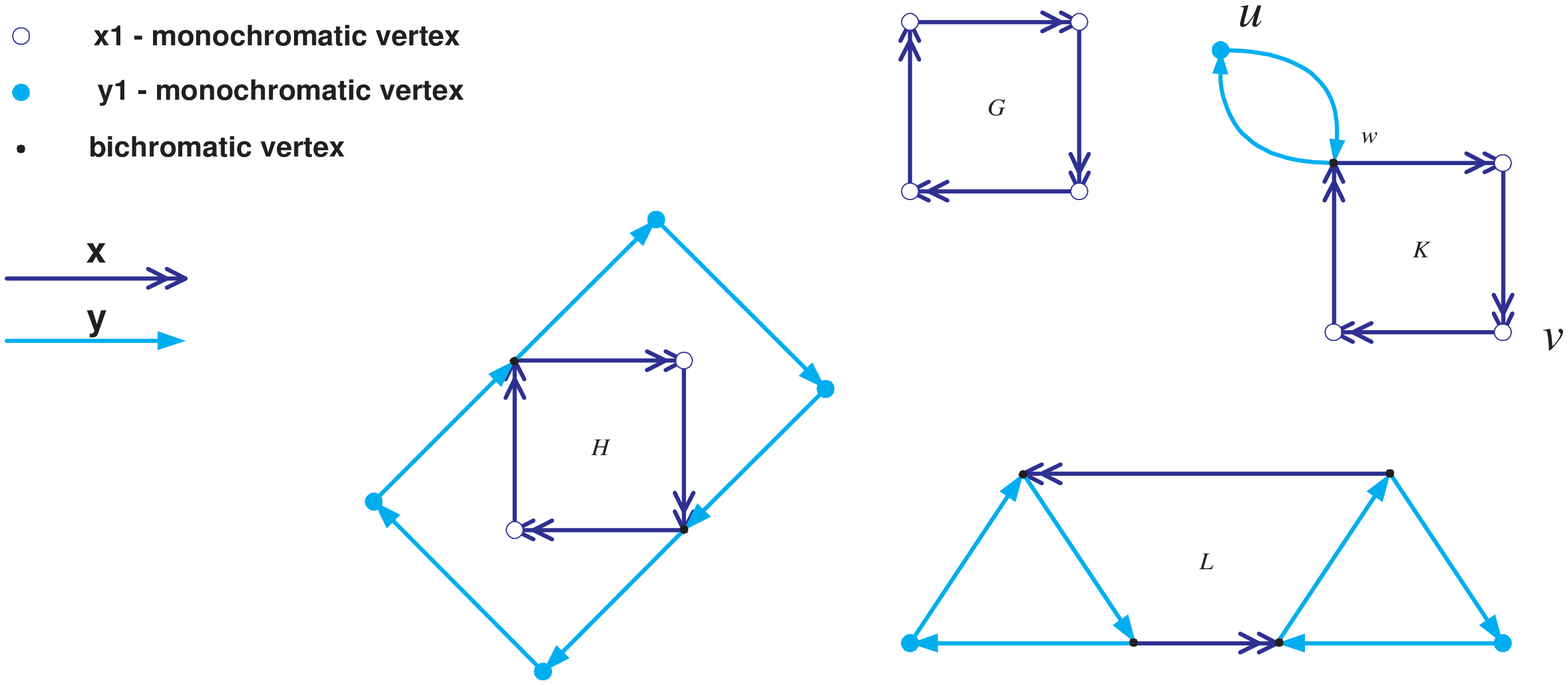}
\caption{ \label{fig:Precovers}}
\end{figure}

\begin{ex}
{\rm
 Let $G=gp\langle x,y | x^4, y^6, x^2=y^3 \rangle=\mathbb{Z}_4 \ast_{\mathbb{Z}_2} \mathbb{Z}_6$.

Recall that $G$ is isomorphic to $SL(2,\mathbb{Z})$ under the
homomorphism
$$x\mapsto \left(
\begin{array}{cc}
0 & 1 \\
-1 & 0
\end{array}
\right), \ y \mapsto \left(
\begin{array} {cc}
0 & -1\\
1 & 1
\end{array}
 \right).$$
The graphs $\Gamma_1$ and $\Gamma_3$ on Figure \ref{fig:Precovers}
are examples of precovers of $G$ with one monochromatic component
and two monochromatic components, respectively.

Though the $\{x\}$-monochromatic component of the graph $\Gamma_2$
is a cover of $\mathbb{Z}_4 $ and the $\{y\}$-monochromatic
component is a cover of $\mathbb{Z}_6$, $\Gamma_2$ is not a
precover of $G$, because it is not a $G$-based graph. Indeed, $v
\cdot (x^2y^{-3})=u$, while $x^2y^{-3}=_G 1$.

The graph $\Gamma_4$ is not a precover of $G$ because its
$\{x\}$-monochromatic components are not covers of  $\mathbb{Z}_4
$. }\e
\end{ex}


A graph $\Gamma$ is \emph{$x$-saturated} at $v \in V(\Gamma)$, if
there exists $e \in E(\Gamma)$ with $\iota(e)=v$ and $lab(e)=x$
($x \in X$). $\Gamma$ is \emph{$X^{\pm}$-saturated} if it is
$x$-saturated for each $x \in X^{\pm}$ at each  $v \in V(\Gamma)$.

\begin{lem}[Lemma 1.5 in \cite{gi_sep}] \label{lemma1.5}
Let $G=gp\langle X|R \rangle$ be a group and let $(\Gamma,v_0)$ be
a graph well-labelled with $X^{\pm}$. Denote $Lab(\Gamma,v_0)=S$.
Then
\begin{itemize}
    \item $\Gamma$ is $G$-based if and only if it can be embedded in $(Cayley(G,S), S~\cdot~1)$,
    \item $\Gamma$ is $G$-based and $X^{\pm}$-saturated if and only if it is isomorphic to  \linebreak[4] $(Cayley(G,S), S
    \cdot~1)$.~
    \footnote{We write $S \cdot 1$ instead of the usual $S1=S$ to distinguish this vertex of $Cayley(G,S)$ as the basepoint of the
    graph.}
\end{itemize}
\end{lem}

\begin{cor} \label{cor:PrecoversSubgrOfCayleyGr}
If $\Gamma$ is a precover of $G$ with $Lab(\Gamma,v_0)=H \leq G$
then  $\Gamma$ is a subgraph of $Cayley(G,H)$.
\end{cor}

Thus a precover of $G$ can be viewed as a part of the
corresponding cover  of $G$, which explains the use of the term
``precovers''.

\begin{remark}[\cite{m-foldings}] \label{remark: morphism of precovers}
{\rm Let $\phi: \Gamma \rightarrow \Delta$ be a morphism of
labelled graphs. If $\Gamma$ is a precover of $G$, then
$\phi(\Gamma)$ is a precover of $G$ as well. }\e
\end{remark}


\subsection*{Precovers are Compatible:}

A graph $\Gamma$  is called \emph{compatible at a bichromatic
vertex} $v$ if for any monochromatic path $p$ in $\Gamma$ such
that $\iota(p)=v$ and $lab(p) \in A$ there exists a monochromatic
path $t$ of a different color in $\Gamma$ such that $\iota(t)=v$,
$\tau(t)=\tau(p)$ and $lab(t)=_G lab(p)$. We say that $\Gamma$ is
\emph{compatible} if it is compatible at all bichromatic vertices.

\begin{ex}
{\rm The graphs $\Gamma_1$ and $\Gamma_3$ on Figure
\ref{fig:Precovers} are compatible. The graph $\Gamma_2$ does not
possess this property because $w \cdot x^{2}=v$, while $w \cdot
y^3=u$. $\Gamma_4$ is not compatible as well.} \e
\end{ex}

\begin{lem} [Lemma 2.12 in \cite{gi_sep}] \label{lemma2.12}
If $\Gamma$ is a compatible graph, then for any  path $p$ in
$\Gamma$ there exists a path $t$ in normal form  such that
$\iota(t)=\iota(p), \ \tau(t)=\tau(p) \ {\rm and} \ lab(t)=_G
lab(p).$
\end{lem}

\begin{remark} [Remark 2.11 in \cite{gi_sep}] \label{remark:
precovers are compatible} {\rm Precovers are compatible. \hfill
$\diamond$}
\end{remark}

The following can be taken as another definition of precovers.

\begin{lem} [Corollary2.13 in \cite{gi_sep}]  \label{corol2.13}
Let $\Gamma$ be a compatible graph. If all  $X_i$-components of
$\Gamma$ are $G_i$-based, $i \in \{1,2\}$, then $\Gamma$ is
$G$-based. In particular, if each $X_i$-component of $\Gamma$ is a
cover of $G_i$, $i \in \{1,2\}$, and $\Gamma$ is compatible, then
$\Gamma$ is a precover of $G$.
\end{lem}


\subsection*{Complexity Issues:}
As were noted in \cite{m-foldings}, the complexity of the
generalized Stallings' algorithm is quadratic in the size of the
input, when we assume that all the information concerning the
finite groups $G_1$, $G_2$, $A$ and the amalgam $G=G_1 \ast_{A}
G_2$ given via $(1.a)$, $(1.b)$ and $(1.c)$ (see
Section~\ref{section:Preliminaries}) is not a part of the input.
We also assume that the Cayley graphs and all the relative Cayley
graphs of the free factors are given for ``free'' as well.

Otherwise, if the group presentations of the free factors $G_1$
and $G_2$, as well as the monomorphisms between the amalgamated
subgroup $A$ and the free factors are a part of the input (the
\emph{uniform version} of the algorithm) then we have to build the
groups $G_1$ and $G_2$, that is to construct their Cayley graphs
and relative Cayley graphs.

Since we assume that the groups $G_1$ and $G_2$ are finite,  the
Todd-Coxeter algorithm and the Knuth Bendix algorithm  are
suitable \cite{l_s, sims, stil} for these purposes. Then the
complexity of the construction depends on the group presentation
of $G_1$ and $G_2$ we have: it could be even exponential in the
size of the presentation \cite{cdhw73}. Therefore the generalized
Stallings algorithm, presented in \cite{m-foldings}, with these
additional constructions could take time exponential in the size
of the input.

Thus each uniform algorithmic problem for $H$ whose solution
involves the construction of the subgroup graph $\Gamma(H)$ may
have an exponential complexity in the size of the input.

The primary goal of the complexity analysis introduced along the
current paper is to estimate our graph theoretical methods. To
this end,  we assume that all the algorithms along the present
paper have the following ``given data''.
\begin{description}
    \item[GIVEN] : Finite groups $G_1$, $G_2$, $A$ and the amalgam
$G=G_1 \ast_{A} G_2$ given via $(1.a)$, $(1.b)$ and $(1.c)$.\\
We assume that the Cayley graphs and all the relative Cayley
graphs of the free factors are given.
\end{description}


\section{Computing Subgroup Presentations}
\label{sec:SubgroupPresentation}

Given a presentation of a group $G$, and a suitable information
about a subgroup $H$ of $G$, the Reidemeister-Schreier method (see
2.2.3 in \cite{m-k-s}) enables one  to compute a presentation for
$H$.

It's a well known fact (see, for instance, \cite{m-k-s}, p.90)
that if $[G:H]<\infty$ then $H$ is finitely generated when $G$ is
finitely generated, and $H$ is finitely presented when $G$ is
finitely presented. Such a finite presentation of $H$ can be
effectively calculated by an application of the
Reidemeister-Schreier method.

However a subgroup can be finitely presented even if its  index is
infinite. For instance, if the group under the consideration is
\emph{coherent}, then all its finitely generated subgroups are
finitely presented. Recently, coherence of some classes of groups
has been investigated in \cite{m_w, p_s}.

Below we introduce a restricted version of the
Reidemeister-Schreier method which allows to compute a finite
presentation for a finitely generated subgroup $H$ of an amalgam
of finite groups $G=G_1 \ast_A G_2$ given by (1.a)-(1.c). This
immediately implies the \emph{coherence} of amalgams of finite
groups.

The suitable information about the subgroup which is needed for an
application of the method  can be read off from its subgroup graph
$\Gamma(H)$ constructed by the generalized Stallings' algorithm.

Let $(\Gamma,v_0)$ be a finite precover of $G$. Let
$H=Lab(\Gamma,v_0)$.

Recall that $H=Lab(\Gamma,v_0)$ is the image of
$lab(Loop(\Gamma,v_0)) \subseteq X^*$ in $G$ under the natural
morphism  $\varphi: X^* {\rightarrow} G$. Note that
$\varphi=\varphi_2 \circ \varphi_1$, where
$$\varphi_1: X^*
\rightarrow FG(X) \ {\rm and} \ \varphi_2: FG(X) \rightarrow G.$$
Let  $\widetilde{H} = \varphi_1 \left(lab(Loop(\Gamma,v_0)
\right)$. Thus  $H=\varphi_2(\widetilde{H})$.
Moreover, $$H=\widetilde{H}  / N=\widetilde{H}  /
\left(\widetilde{H} \cap N \right),$$ where $N$ is the normal
closure of $R$ in $FG(X)$ (see \cite{l_s, m-k-s}). We put
$F=FG(X)$.

Let $T$ be a fixed spanning tree  of $\Gamma$. For all $v \in
V(\Gamma)$, we consider $t_v$ to be the unique freely reduced path
in $T$ from the basepoint $v_0$ to the vertex $v$.

For each $e \in E(\Gamma)$ we consider $t(e)=t_{\iota(e)}e
\overline{t_{\tau(e)}}$. Thus if $e \in E(T)$ then $t(e)$ can be
freely reduced to an empty path, that is $lab(t(e))=_{F} 1$.

Let $E^+$ be the set of positively oriented edges of $\Gamma$. Let
\begin{equation} \label{eq:Def_of X_H}
X_H=\{lab(t(e)) \; | \; e \in E^+ \setminus E(T) \},
\end{equation}
$$Q_v=\{ q \subseteq \Gamma  \; | \; \iota(q)=\tau(q)=v, \;
lab(q) \equiv r \in R \},$$
\begin{equation} \label{eq:DefI of R H}
R_H=\left\{lab\left(\phi\left(t_vq\overline{t_v}\right)\right) \;
|  v \in V(\Gamma), \; t_v \subseteq T, \; q \in Q_v,\right\},
\end{equation}
where   $\phi$ is a function from the set of freely reduced paths
in $\Gamma$ into $Loop(\Gamma,v_0)$ defined as follows.
$$\phi(p)=t(e_1)t(e_2) \cdots t(e_n), \ {\rm where} \
 p=e_1e_2 \cdots e_n \subseteq \Gamma.$$   Thus the
path $\phi(p)$ is closed at $v_0$ in $\Gamma$ and
$$lab(\phi(p)) \equiv lab(t(e_1))lab(t(e_2)) \cdots lab(t(e_n)).$$
Moreover, if the path $p$ is closed at $v_0$ in $\Gamma$ then the
path $\phi(p)$ is \emph{freely equivalent} to $p$, that is
$\phi(p)$ can be transformed to the path $p$ by a series of free
reductions. Thus $lab(\phi(p))=_{F} lab(p)$.

The function $\phi$ induces a partial function $\phi'$ from
$FG(X)$ into $FG(X_H)$ such that $\phi'(w) = lab(\phi(p))$, where
$p$ is a path in $\Gamma$ with $lab(p) \equiv w$. Thus another
definition of $R_H$ takes the following form
\begin{equation} \label{eq:DefII of R H}
R_H=\left\{\phi'\left(lab(t_vq\overline{t_v})\right)
\; |  \; v \in V(\Gamma), \; t_v \subseteq T, \; q \in
Q_v,\right\}.
\end{equation}

\begin{remark}
{\rm Note that the system of coset representatives $ \{lab (t_v)
\; | \; v \in V(\Gamma) \}$ is a subset of the \emph{Schreier
transversal} of $\widetilde{H}$ in $FG(X)$ (\cite{stal}).  }\e
\end{remark}


\begin{thm} \label{thm:NewReidmeisterSchreier}

With the above notation,   $H=gp \langle X_H \; | R_H \rangle$.
\end{thm}
\begin{proof}
As is well known (\cite{kap-m, mar_meak, stal}),
$\widetilde{H}=FG(X_H)$. Therefore $H=\langle X_H \rangle$.

To complete the proof it remains to show that   the normal closure
$N_H$ of $R_H$ in $FG(X_H)=\widetilde{H}$ is equal to
$\widetilde{H}  \cap N$.

Let $v \in V(\Gamma)$ such that $Q_v \neq \emptyset$. Let $q \in
Q_v$. Therefore $\phi(t_vq\overline{t_v})$ is freely equivalent to
the path $t_vq\overline{t_v}$. Thus
$$lab(\phi(t_vq\overline{t_v}))=_{F}lab(t_vq\overline{t_v}) \equiv lab(t_v)lab(q)lab(t_v)^{-1}
\in N.$$ On the other hand, the path $t_vq\overline{t_v}$ is
closed at $v_0$, hence $lab(t_vq\overline{t_v}) \in
\widetilde{H}$. Thus $lab(\phi(t_vq\overline{t_v})) \in
\widetilde{H} \cap N$. Therefore $R_H \subseteq \widetilde{H} \cap
N$.

For all $y \in \widetilde{H}$, there exist a closed path $s \in
\Gamma$ starting at $v_0$ with $lab(s)=_F \widetilde{H}$. By the
definition of $R_H$, for all $r \in R_H$ there exist a path
$t_vq\overline{t_v} \subseteq \Gamma$ closed at $v_0$ such that
$lab(t_vq\overline{t_v})=_F r$. Hence the path
$s(t_vq\overline{t_v})\overline{s}$ is closed at $v_0$ in
$\Gamma$. Thus $lab(s(t_vq\overline{t_v})\overline{s})=_F yry^{-1}
\in \widetilde{H}$. Moreover,
$lab(s(t_vq\overline{t_v})\overline{s}) \equiv
lab(st_v)lab(q)(lab(st_v))^{-1} \in N$. Therefore $N_H \subseteq
\widetilde{H} \cap N$.

{ \ }

Assume now that $w \in \widetilde{H} \cap N$.
%
%
Since $w \in \widetilde{H}$, there exists  a freely reduced path
$p$ in $\Gamma$ closed at $v_0$ with $lab(p)=_{F} w$ (\cite{kap-m,
mar_meak}). Let $p=p_1 \cdots p_k$ be its decomposition   into
maximal monochromatic paths $p_i$ with $lab(p_i) \equiv w_i \in
G_{l_i}$ ($1 \leq i \leq k$ and $l_i \in \{1,2\}$).

Since $w \in N$, $w=_G 1$. Therefore, by the Normal Form Theorem
for  free products with amalgamation (IV.2.6 in \cite{l_s}), there
exists $1 \leq i \leq k$ such that  $w_i   \in A \cap G_{l_i}$.
The proof is by induction on the number $k$ of the maximal
monochromatic subpaths of the path $p$. Without loss of
generality,  simplifying the notation, we let $l_i=1$.


Assume first that $w_i=_{G_1} 1$. Since $\Gamma$ is $G$-based, the
subpath $p_i$ is closed at $\iota(p_i)=\tau(p_i)$.

Let $v_j =\tau(p_j)$ and let $t_j=t_{v_j} \subseteq T$ ($1 \leq j
\leq k$). Thus $v_{i-1}=v_i$ and $t_{i-1}=t_i$. Let $t=p_1 \cdots
p_{i-1}$. See Figure~\ref{fig:ProofReidmSchreier} (a).


\begin{figure}[!h]
\begin{center}
\psfrag{v0 }[][]{$v_0$}
\psfrag{v1 }[][]{$v_1$}
\psfrag{vk }[][]{$v_{k-1}$}
\psfrag{v }[][]{$v_{i-1}$}
\psfrag{p1 }[][]{$p_1$}
\psfrag{pk }[][]{$p_k$}
\psfrag{px }[][]{$p_{i-1}$}
\psfrag{py }[][]{$p_{i+1}$}
\psfrag{pi }[][]{$p_i$}
\psfrag{pj }[][]{$p_i'$}
\psfrag{ti }[][]{$t_{i-1}$}
\psfrag{t }[][]{$t$}

\includegraphics[width=\textwidth]{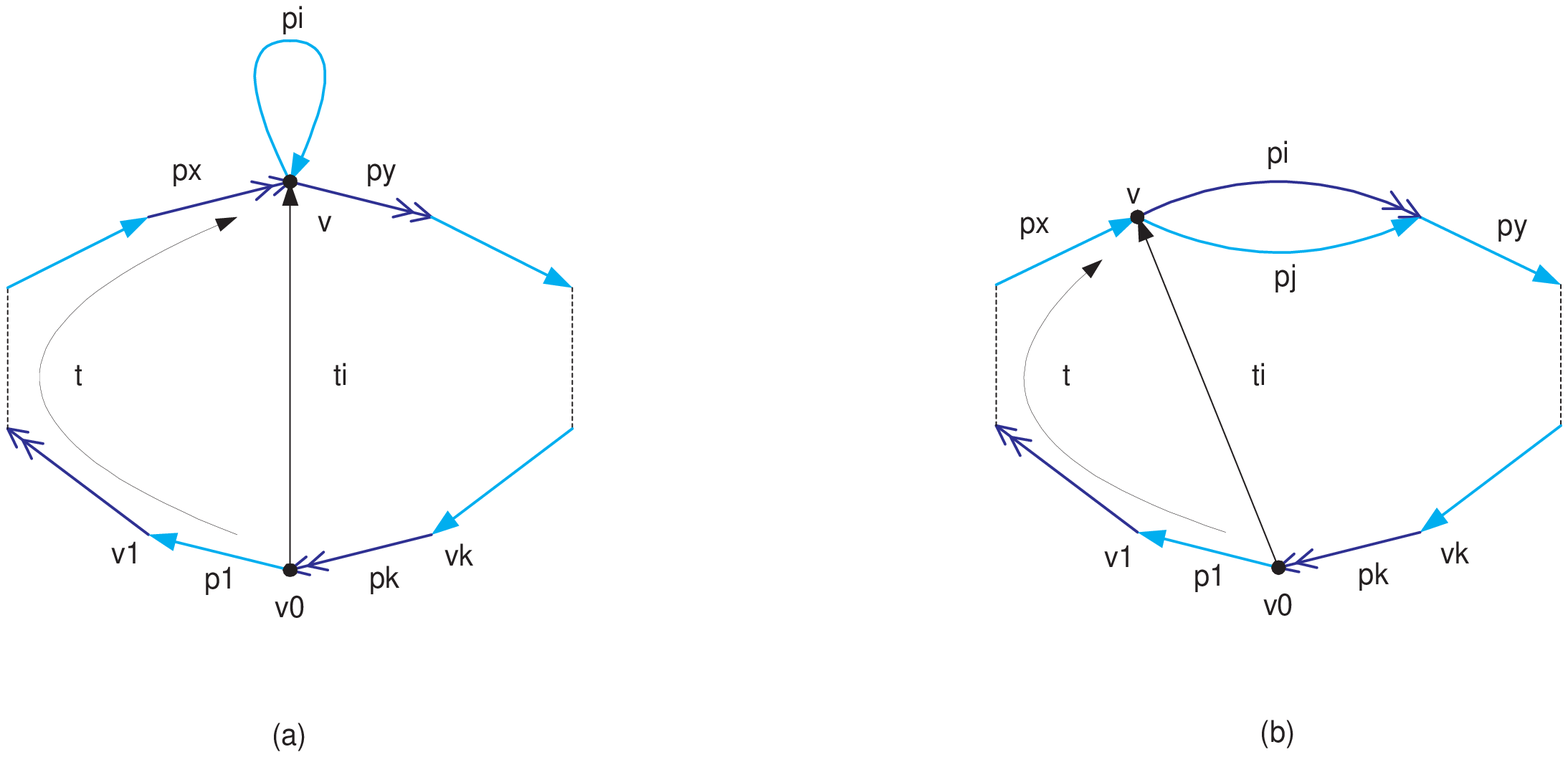}
\caption{ \label{fig:ProofReidmSchreier}}
\end{center}
\end{figure}


Hence the path $p$ can be obtained by free reductions from the
following path
$$ \left( \: (t\overline{ t}_{i-1}) (t_{i-1} p_i
\overline{t_{i-1}}) (t_{i-1}\overline{t}) \: \right)\left(tp_{i+1}
\cdots p_k \right).$$

Thus $lab(t\overline{ t}_{i-1}) \in \widetilde{H}$, and the number
of the maximal monochromatic subpaths of the path
$$tp_{i+1} \cdots p_k=p_1 \cdots p_{i-2}(p_{i-1}p_{i+1})p_{i+2}
\cdots p_k$$ is  $k-2$. Therefore, by the inductive assumption,
$lab (tp_{i+1} \cdots p_k) \in N_H$. To get the desired conclusion
it remains to show that $lab(t_{i-1} p_i \overline{t_{i-1}}) \in
N_H$.

Since $w_i=_{G_1} 1$, we have $lab(p_i)  \in N_{1}$, where $N_{1}$
is the normal closure of $R_{1}$ in $F_1=FG(X_{1})$. Therefore
$lab(p_i) =_{F_1} (z_{1}s_{1}z_{1}^{-1}) \cdots
(z_{m}s_{m}z_{m}^{-1})$, where $z_{j} \in F_1$ and $s_{j} \in
R_{1}$ ($1 \leq j \leq m$).

Let  $C_i$ be the $X_{1}$-monochromatic component of $\Gamma$ such
that $p_i \subseteq C_i$. Since $\Gamma$ is a precover of $G$, $C$
is $X_{1}^{\pm}$-saturated. Hence $p_i$ is a free reduction of the
path $p_i' \subseteq C$ such that $lab(p_i') \equiv
(z_{1}s_{1}z_{1}^{-1}) \cdots (z_{m}s_{m}z_{m}^{-1})$. Since
$\Gamma$ is $G$-based, the subpaths of $p_i'$ labelled by $s_j$
($1 \leq j \leq m$) are closed. Therefore $p_i'$ has the following
decomposition
$$p_i'=(c_{1}q_{1}\overline{c_{1}}) \cdots
(c_{m}q_{m}\overline{c_{m}}),$$ where $lab(c_{j}) \equiv z_{j}$
and $lab(q_{j}) \equiv s_{j}$ ($1 \leq j \leq m$). Thus the path
$t_{i-1} p_i \overline{t_{i-1}}$ can be obtained by free
reductions from the path
$$\left(t_{i-1}(c_{1}q_{1}\overline{c_{1}})\overline{t_{i-1}}\right) \cdots
\left(t_{i-1}(c_{m}q_{m}\overline{c_{m}})\overline{t_{i-1}}\right).$$
For each $1\leq j \leq m$, the path
$t_{i-1}(c_{j}q_{j}\overline{c_{j}})\overline{t_{i-1}}$ is a free
reduction of the path
$$(t_{i-1} c_{1}
\overline{t_{\tau(c_{j})}})(t_{\tau(c_{j})}q_{j}\overline{t_{\tau(c_{j})}})(\overline{t_{i-1}
c_{j} \overline{t_{\tau(c_{j})}}}).$$

Since $lab(t_{i-1} c_{j} \overline{t_{\tau(c_{j})}}) \in
\widetilde{H}$ and
$lab(\phi(t_{\tau(c_{j})}q_{j}\overline{t_{\tau(c_{j})}})) \in
R_H$ we conclude that
$lab(t_{i-1}(c_{j}q_{j}\overline{c_{j}})\overline{t_{i-1}}) \in
N_H$ ($1 \leq j \leq m$). Therefore $lab(t_{i-1} p_i
\overline{t_{i-1}}) \in N_H$. We are done.

{ \ } \\


Assume now that   $1 \neq_G w_i \in A \cap G_{1}$.

Since for all $1 \leq j \leq k$ the vertices $v_j=\tau(p_j)$ are
bichromatic, and because the graph $\Gamma$ is compatible, there
exists a $X_2$-monochromatic path $p_i'$ in $\Gamma$ such that
$\iota(p_i')=\iota(p_i)$, $\tau(p_i')=\tau(p_i)$ and $lab(p_i')=_G
lab(p_i)$. See Figure~\ref{fig:ProofReidmSchreier} (b).
%
%
%
Hence the path $p$ can be obtained by free reductions from the
following path
$$ \left( \: (t\overline{ t}_{i-1}) (t_{i-1} p_i\overline{p_i'}
\overline{t_{i-1}}) (t_{i-1}\overline{t}) \: \right)
\left(tp'_ip_{i+1} \cdots p_k \right).$$

Thus $lab(t\overline{ t}_{i-1}) \in \widetilde{H}$, and the number
of the maximal monochromatic subpaths of the path
$$tp'_ip_{i+1} \cdots p_k=p_1 \cdots p_{i-2}(p_{i-1}p'_ip_{i+1})p_{i+2}
\cdots p_k$$ is  $k-2$. Therefore, by the inductive assumption,
$lab (tp'_ip_{i+1} \cdots p_k) \in N_H$. To get the desired
conclusion it remains to show that $lab(t_{i-1}
(p_i\overline{p'_i}) \overline{t_{i-1}}) \in N_H$.

Let $lab(p_i) =_{G_1} a_1  \cdots a_m $, where $a_j$ are
generators of $A \cap G_1$. Let $b_j$ be corresponding generators
of $A \cap G_2$ such that $a_j=_G b_j$ and $a_j{b_j}^{-1} \in R$
($1 \leq j \leq m$). Note that
$$ (a_1  \cdots a_m)(b_1  \cdots b_m)^{-1}   =_{F} $$
$$ =_F  \left(a_1{b_1}^{-1}\right)\left(b_1 \left( a_2{b_2 }^{-1}
\right) b_1 ^{-1}\right) \cdots
\left(b_1  \cdots b_{m-1} \left(a_m b_m ^{-1} \right) b_{m-1}
^{-1} \cdots b_1^{-1}\right).$$

Since monochromatic components of $\Gamma$ are
$X_i^{\pm}$-saturated ($i \in \{1,2\}$), and  because $\iota(p_i)
\in VB(\Gamma)$, there exist paths $\gamma_1$ and $\delta_1$ such
that $\iota(\gamma_1)=\iota(p_i)=\iota(\delta_1)$ and $
lab(\gamma_1) \equiv a_1$, $lab(\delta_1) \equiv b_1$. Since
$\Gamma$ is compatible, $\tau(\gamma_1)=\tau(\delta_1) \in
VB(\Gamma)$. Thus there exist paths $\gamma_2$ and $\delta_2$ such
that $\iota(\gamma_2)=\tau(\gamma_1)=\iota(\delta_2)$ and $
lab(\gamma_2) \equiv a_2$, $lab(\delta_2) \equiv b_2$. Since
$\Gamma$ is compatible, $\tau(\gamma_2)=\tau(\delta_2) \in
VB(\Gamma)$.

Continuing  in this manner one can construct such paths
$\gamma_j$, $\delta_j$ for all $1 \leq j \leq m$. Thus $p_i$ and
$p'_i$  are free reductions of the paths $\gamma_1 \cdots
\gamma_m$ and $\delta_1 \cdots \delta_m$, respectively.
Hence the path $p_i\overline{p_i'}$ can be obtained by free
reductions from the path
$$  \left(\gamma_1{\overline{\delta_1}} \right)\left(\delta_1 \left( \gamma_2{\overline{\delta_2}
} \right) \overline{\delta_1}  \right) \cdots
\left(\delta_1  \cdots \delta_{m-1} \left(\gamma_m
\overline{\delta_m } \right)\overline{ \delta_{m-1}}  \cdots
\overline{\delta_1} \right).$$
Therefore the path $t_{i-1}(p_i\overline{p_i'})\overline{t}_{i-1}$
is a free reduction of
$$  \left(t_{i-1}(\gamma_1{\overline{\delta_1}}) \overline{t}_{i-1} \right)
  \cdots
\left(t_{i-1}(\delta_1  \cdots \delta_{m-1} \left(\gamma_m
\overline{\delta_m } \right)\overline{ \delta_{m-1}}  \cdots
\overline{\delta_1}) \overline{t}_{i-1}\right).$$

For each $1 \leq j \leq m-1$, the path $t_{i-1}(\delta_1  \cdots
\delta_{j-1} \left(\gamma_j \overline{\delta_j } \right)\overline{
\delta_{j-1}}  \cdots \overline{\delta_1}) \overline{t}_{i-1}$ is
a free reduction of the path
$$\left(t_{i-1}\delta_1 \cdots \delta_{j-1}
\overline{t_{\iota(\gamma_j)}}\right)
\left(t_{\iota(\gamma_j)}(\gamma_j{\overline{\delta_j} })
\overline{t_{\iota(\gamma_j)}} \right) \left(t_{\iota(\gamma_j)}
\overline{\delta_{j-1}} \cdots \overline{\delta_1}
\overline{t_{i-1}}\right).$$

Since  $ lab\left(\phi\left(t_{\iota(\gamma_j)}(\gamma_j
\overline{\delta_j}) \overline{t_{\iota(\gamma_j)}}\right)\right)
\in R_H$   and
$lab(t_{i-1}\delta_1 \cdots \delta_{j-1}
\overline{t_{\iota(\gamma_j)}}) \in \widetilde{H}$, we conclude
that  for each $1 \leq j \leq m-1$
$$lab(t_{i-1}(\delta_1 \cdots
\delta_{j-1} \left(\gamma_j \overline{\delta_j } \right)\overline{
\delta_{j-1}}  \cdots \overline{\delta_1}) \overline{t}_{i-1}) \in
N_H.$$
Therefore $lab(t_{i-1} (p_i\overline{p_i'}) \overline{t_{i-1}})
\in N_H$. We are done.

\end{proof}

\begin{cor} \label{cor:ReidShcrPrecovers}
Let $(\Gamma,v_0)$ be a finite precover of $G$.
Then there exists an algorithm which computes a subgroup of $G$
determined by $(\Gamma,v_0)$, that is computes a finite group
presentation of $H=Lab(\Gamma,v_0)$.
\end{cor}
\begin{proof}
We  compute the sets $X_H$ and $R_H$ according to their
definitions. These sets are finite, because the graph $\Gamma$ is
finite.
By Theorem \ref{thm:NewReidmeisterSchreier}, $H= gp\langle X_H \;
| \; R_H \rangle$.

\end{proof}


\begin{cor} \label{cor:NewReidmShcreier}
Let $h_1, \ldots h_n \in G$.
Then there exists an algorithm which computes a finite group
presentation of the subgroup $H=\langle h_1, \ldots, h_n \rangle$
in $G$ (not necessary with respect to $\{h_1, \cdots, h_n\}$) .
\end{cor}
\begin{proof} We first construct the graph $(\Gamma(H),v_0)$, using the
generalized Stallings' folding  algorithm. By Theorem~\ref{thm:
properties of subgroup graphs} (2), this graph is a finite
precover of $G$.
Now we proceed according to Corollary~\ref{cor:ReidShcrPrecovers}.

\end{proof}

\begin{cor}
Amalgams of finite groups are coherent.
\end{cor}

\begin{remark}
{\rm As is well known, the Reidemeister-Schreier method yields a
presentation of a subgroup $H$ which is usually not in a useful
form. Namely, some of the generators are redundant and can be
eliminated, while some of the relators can be simplified. In order
to improve (to simplify) this presentation, one can apply the
Tietze transformation. An efficient version of such a
simplification procedure was developed in \cite{havas, hkrr}.
 }\e
\end{remark}

\begin{figure}[!htb]
\psfrag{v0 }[][]{$v_0$} \psfrag{x }[][]{$x$}
 \psfrag{y }[][]{$y$}

\includegraphics[width=0.5\textwidth]{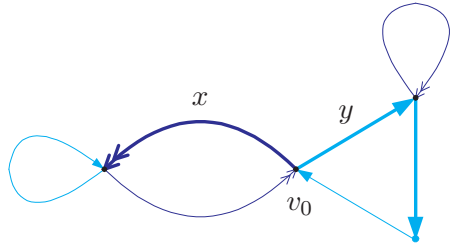}
\caption {{\footnotesize The bold edges of the  graph $\Gamma(H)$
correspond to a spanning tree $T$} \label{fig:ExReidmSchreier}}
\end{figure}

\begin{ex} \label{ex:ReidmShreier}
{\rm Let $G=gp\langle x,y | x^4, y^6, x^2(y^3)^{-1}
\rangle=\mathbb{Z}_4 \ast_{\mathbb{Z}_2} \mathbb{Z}_6$.

Recall that $G$ is isomorphic to $SL(2,\mathbb{Z})$ under the
homomorphism
$$x\mapsto \left(
\begin{array}{cc}
0 & 1 \\
-1 & 0
\end{array}
\right), \ y \mapsto \left(
\begin{array} {cc}
0 & -1\\
1 & 1
\end{array}
 \right).$$

Let $H  =\langle xyx^{-1}, yxy^{-1} \rangle$ be  a subgroup of
$G$. The subgroup graph $\Gamma(H)$ constructed by the generalized
Stallings' folding  algorithm is presented on
Figure~\ref{fig:ExReidmSchreier}.

We apply to $\Gamma(H)$ the algorithm described along with the
proof of Corollary~\ref{cor:NewReidmShcreier}.

We first compute $X_H$ according to (\ref{eq:Def_of X_H}):
$$h_1=xyx^{-1}, \ h_2=x^2, \ h_3=yxy^{-1}, \ h_4=y^3.$$

The computation of $R_H$  according to (\ref{eq:DefII of R H})
consists of the following steps.

$$\phi'(x^4)=(h_2)^2, \ \phi'(y^6)=(h_4)^2, \
\phi'(x^2(y^3)^{-1})=h_2(h_4)^{-1}.$$

$$\phi'(x(x^4)x^{-1})=(h_2)^2, \ \phi'(x(y^6)x^{-1})=(h_1)^6, \
\phi'(x(x^2(y^3)^{-1})x^{-1})=h_2(h_1)^{-3}.$$

$$\phi'(y(x^4)y^{-1})=(h_3)^4, \ \phi'(y(y^6)y^{-1})=(h_4)^2, \
\phi'(y(x^2(y^3)^{-1})y^{-1})=h_3^2(h_4)^{-1}.$$

Therefore $H  =gp\langle h_1, h_3 \; | \; h_1^6, h_3^4,
h_1^3=h_3^2 \rangle$.

}\e
\end{ex}

\subsection*{Complexity.}

Let $m$ be the sum of the lengths of the words $h_1, \ldots h_n$.
By Theorem~\ref{thm: properties of subgroup graphs} $(4)$, the
generalized Stallings' algorithm  computes $(\Gamma(H),v_0)$ in
time $O(m^2)$.

The construction of $X_H$, which is a free basis of
$\widetilde{H}= \varphi_1 \left(lab(Loop(\Gamma(H),v_0) \right)$,
takes $O(|E(\Gamma(H))|^2)$, by \cite{b-m-m-w}. Since, by
Theorem~\ref{thm: properties of subgroup graphs} $(4)$,
$|E(\Gamma(H))|$ is proportional to $m$, the computation of $X_H$
takes $O(m^2)$.

To construct the set $Q_H$ we try to read each one of the defining
relators of $G$ at each one of the vertices of the graph
$\Gamma(H)$. It takes at most $$|R| \cdot |V(\Gamma(H)| \cdot
\left(\sum_{v \in V(\Gamma(H))} deg(v) \right).$$ Since $\sum_{v
\in V(\Gamma(H))} deg(v)=2 |E(\Gamma(H)|$  and because, by our
assumption, the presentation of $G$ is given and it is not a part
of the input, the computation of the set $Q_H$ takes
$O(|V(\Gamma(H)| \cdot |E(\Gamma(H)|)$. Since, by
Theorem~\ref{thm: properties of subgroup graphs} (4),
$|V(\Gamma(H)|=O(m)$, it takes $O(m^2)$.

The rewriting process which yield the set of relators $R_H$ takes
at most $|V(\Gamma(H)| \cdot \left(\sum_{r \in R} |r| \right)$
which is $O(|V(\Gamma(H)|)$.

Thus the complexity of the restricted Reidemeister-Schreier
process given by Corollary~\ref{cor:NewReidmShcreier} is $O(m^2)$.


\section{The Freeness Problem}
\label{sec:FreenessProblem}

A  freeness of subgroups is one of the fundamental questions of
combinatorial and geometric group theory. The classical results in
this issue include the Nielsen-Schreier subgroup theorem for free
groups, the corollary of Kurosh subgroup theorem and the
Freiheitssatz of Magnus.

Namely, subgroups of free groups are free (I.3.8, \cite{l_s}). A
subgroup of a free product which has a trivial intersection with
all conjugates of the factors is free (\cite{l_s}, p.120). A
subgroup $H$ of an one-relator group $G=gp\langle X \; | \; r=1
\rangle$, where $r$ is \emph{cyclically freely reduced}, is free
if $H$ is generated by a subset of $X$ which  omits a generator
occurring in $r$ (II.5.1, \cite{l_s}).

Results concerning amalgamated free products follow from the
Neumann's subgroup theorem.
\begin{thm}[H.Neumann, IV.6.6 \cite{l_s}]
Let $G=G_1 \ast_A G_2$ be a non-trivial free product with
amalgamation. Let $H$ be a finitely generated subgroup of $G$ such
that all conjugates of $H$ intersect $A$ trivially.

Then $H=F \ast ( \ast_j \; g_jH_jg_j^{-1})$, where  $F$ is a free
group and each $H_j$ is the intersection of a subgroup of $H$ with
a conjugate of a factor of $G$.
\end{thm}

\begin{cor}[IV.6.7 \cite{l_s}] \label{cor: Newmann's thm}
Let $G=G_1 \ast_A G_2$ be a non-trivial free product with
amalgamation. If $H$ is a finitely generated subgroup of $G$ which
has trivial intersection with all conjugates of the factors, $G_1$
and $G_2$, of $G$,  then $H$ is free.
\end{cor}

It turns out (Lemma \ref{free=each component is a Cayley(G_i)})
that the triviality of the intersections between $H$ and
conjugates of the factors, $G_1$ and $G_2$, of $G$ can be detected
from the subgroup graph $\Gamma(H)$ constructed by the generalized
Stallings' folding  algorithm, when $G=G_1 \ast_A G_2$ is an
amalgam of finite groups. Therefore, by Corollary \ref{cor:
Newmann's thm}, the freeness of $H$ is decidable via its subgroup
graph.

We consider the \emph{freeness problem}  to be one which asks to
verify if a  subgroup of a  given group $G$ is free. Clearly, the
freeness problem is solvable in amalgams of finite groups.

Below we introduce a polynomial time algorithm (Corollary
\ref{freeness in amalgams of finite grp}) that employs subgroup
graphs constructed by the generalized Stallings' algorithm to
solve the freeness problem.  A complexity analysis of the
algorithm is given at the end of the section.

\begin{lem} \label{free=each component is a Cayley(G_i)}
Let $H$ be a finitely generated subgroup of an amalgam of finite
groups $G=G_1 \ast_A G_2$.

Then $H$ has a trivial intersection with all conjugates of the
factors  of $G$ if and only if each $X_i$-monochromatic component
$C$ of $\Gamma(H)$ is isomorphic to $Cayley(G_i)$, for all $i \in
\{1,2\}$.
Equivalently, by Lemma~\ref{lemma1.5}, if and only if
$Lab(C,v)=\{1\}$  for each $X_i$-monochromatic component $C$ of
$\Gamma(H)$ ($v \in V(C)$).
\end{lem}
\begin{proof}  Assume first
that there exists a $X_i$-monochromatic component $C$ of
$\Gamma(H)$ ($i \in \{1,2\}$) which is not isomorphic to
$Cayley(G_i)$. Thus, by Lemma~\ref{lemma1.5}, $(C,\vartheta)$ is
not isomorphic to $Cayley(G_i, S, S \cdot 1)$, where is $\vartheta
\in V(C)$ and $\{1\} \neq S \leq G_i$.

Let $1 \neq_G w \in S$. Then there exists a path $q$ in $C$ closed
at $\vartheta$ such that $lab(q) \equiv w$.
Let $p$ be an approach path in $\Gamma(H)$ from $\iota(p)=v_0$ to
$ \tau(p)=\vartheta$. Let $u \equiv lab(p)$.

The path $pq\overline{p}$ is closed at $v_0$ in $\Gamma(H)$. Hence
$lab(pq\overline{p})  \in H$. Therefore $$lab(pq\overline{p})=_G
uwu^{-1} \in H \cap uLab(C,\vartheta)u^{-1}= H \cap uSu^{-1}.$$
Since $w \neq_G 1$, we have $uwu^{-1} \neq_G 1$ and hence $H \cap
uSu^{-1} \neq \{1\}$.

Assume now that there exists $\{1\} \neq S \leq G_i$ ($i \in
\{1,2\}$) such that  $H \cap uSu^{-1} \neq \{1\}$, where $u \in
G$. Let $1 \neq_G h \in H \cap uSu^{-1}$. Thus $h=_G ugu^{-1}$,
where $1 \neq_G g \in S$. Without loss of generality we  can
assume that the words $u$ and $g$ are normal.

If the word $ugu^{-1}$ is in normal form, then there exist a path
$p$ in $\Gamma(H)$ closed at $v_0$ such that $lab(p) \equiv
ugu^{-1}$. Thus there is a decomposition $p=p_1p_2\overline{p_1}$
(because $\Gamma(H)$ is $G$-based, so it is a well-labelled
graph), where $lab(p_1) \equiv u$ and $lab(p_2) \equiv g$. Let $C$
be a $X_i$-monochromatic component of $\Gamma(H)$ such that
$p_2\subseteq C$ and let $v=\tau(p_1)$. Hence $g \equiv lab(p_2)
\in Lab(C, v) \leq G_i$. Thus $Lab(C,v)\neq \{1\}$. Equivalently,
by Lemma~\ref{lemma1.5}, $C$ is not isomorphic to $Cayley(G_i)$.

Assume now that the word $ugu^{-1}$ is not in normal form. Let
$(u_1, \ldots, u_k)$ be a normal decomposition of $u$. Since $g
\in G_i$, its normal decomposition is $(g)$. Hence the normal
decomposition of $ugu^{-1}$ has the form
$$(u_1, \ldots, u_{j-1},w, u_{j-1}^{-1}, \ldots , u_1^{-1}),$$
where $w=_G u_j \ldots u_kgu_k^{-1} \ldots u_{j}^{-1}  \in G_l
\setminus A$ and $u_{j-1} \in G_m \setminus A$ ($1 \leq l \neq m
 \leq 2$).

 Let $u' \equiv u_1 \ldots u_{j-1}$. Then $h=_G u'
w(u')^{-1}$, while the word $u' w(u')^{-1}$ is in normal form and
$w \in G_l$, $l\in\{1,2\}$. Hence, by arguments similar to those
used in the previous case, we are done.

\end{proof}


\begin{thm} \label{cor: h is free}
Let $H$ be a finitely generated subgroup of an amalgam of finite
groups $G=G_1 \ast_A G_2$.

Then $H$ is  free if and only if each $X_i$-monochromatic
component of $\Gamma(H)$ is isomorphic to $Cayley(G_i)$, for all
$i \in \{1,2\}$.
\end{thm}
\begin{proof}
The statement follows immediately  from Corollary \ref{cor:
Newmann's thm} and Lemma \ref{free=each component is a
Cayley(G_i)}.
\end{proof}


Combining Lemma \ref{free=each component is a Cayley(G_i)} with
the Torsion Theorem for amalgamated free products we get Corollary
\ref{cor:ConditionsTorsionFree}.

\begin{thm}[Torsion Theorem, IV.2.7, \cite{l_s}] \label{torsion theorem}
Every element of finite order in $G=G_1 \ast_A G_2$ is a conjugate
of an element of finite order in $G_1$ or $G_2$.
\end{thm}

\begin{cor} \label{cor:ConditionsTorsionFree}
Let $H$ be a finitely generated subgroup of an amalgam of finite
groups $G=G_1 \ast_A G_2$.

Then $H$ is torsion free if and only if each $X_i$-monochromatic
component of $\Gamma(H)$ is isomorphic to $Cayley(G_i)$, for all
$i \in \{1,2\}$.
\end{cor}


\begin{cor} \label{freeness in amalgams of finite grp}
Let $  h_1, \ldots, h_k \in G.$ Then there exists an algorithm
which  decides whether or not the subgroup  $H=\langle h_1,
\ldots, h_k \rangle$ is a free subgroup of $G$.
\end{cor}
\begin{proof} We first construct the graph $\Gamma(H)$, using the
generalized Stallings' folding  algorithm.

 Now, for each $X_i$-monochromatic
component $C$ of $\Gamma(H)$ we verify if $C$ is isomorphic to
$Cayley(G_i)$  ($i \in \{1,2\}$). It can be easily done by
checking the number of vertices of $C$: $|V(C)|=|G_i|$ if and only
if $C$ is isomorphic to $Cayley(G_i)$.

By Theorem \ref{cor: h is free}, $H$ is free if and only if  each
monochromatic component  of $\Gamma(H)$ is isomorphic to the
Cayley graph of an appropriate factor of $G$.

\end{proof}

\begin{remark}
{\rm If $H$ is free then its free basis can be computed using the
restricted Reidemeister-Schreier procedure
(Corollary~\ref{cor:NewReidmShcreier}) followed by a
simplification process based on Tietze transformation. For an
effective version of a simplification procedure when redundant
generators are eliminated consequently using a substring search
technique see \cite{havas, hkrr}. }\e
\end{remark}

%

\begin{ex} \label{ex:Freeness}
{\rm Let $G=gp\langle x,y | x^4, y^6, x^2=y^3 \rangle=\mathbb{Z}_4
\ast_{\mathbb{Z}_2} \mathbb{Z}_6$.

Let $H_1$ and $H_2$ be  finitely generated subgroups of $G$ such
that
$$H_1=\langle xy \rangle \ {\rm and} \ H_2=\langle xy^2, yxyx \rangle.$$

The graphs $\Gamma(H_1)$ and $\Gamma(H_2)$ on Figure \ref{fig:
ExOfFI} are the subgroup graphs of $H_1$ and $H_2$, respectively,
constructed by the generalized Stallings' folding  algorithm. See
Example~\ref{example: graphconstruction} from Appendix  for the
detailed construction of these graphs.

Applying the above algorithm to the graphs $\Gamma(H_1)$ and
$\Gamma(H_2)$, we conclude that $H_2$ is not free, while
$H_1=FG(\{xy\})$. } \e
\end{ex}

\begin{figure}[!htb]
\psfrag{A }[][]{$\Gamma(H_1)$} \psfrag{B }[][]{$\Gamma(H_2)$}

\includegraphics[width=\textwidth]{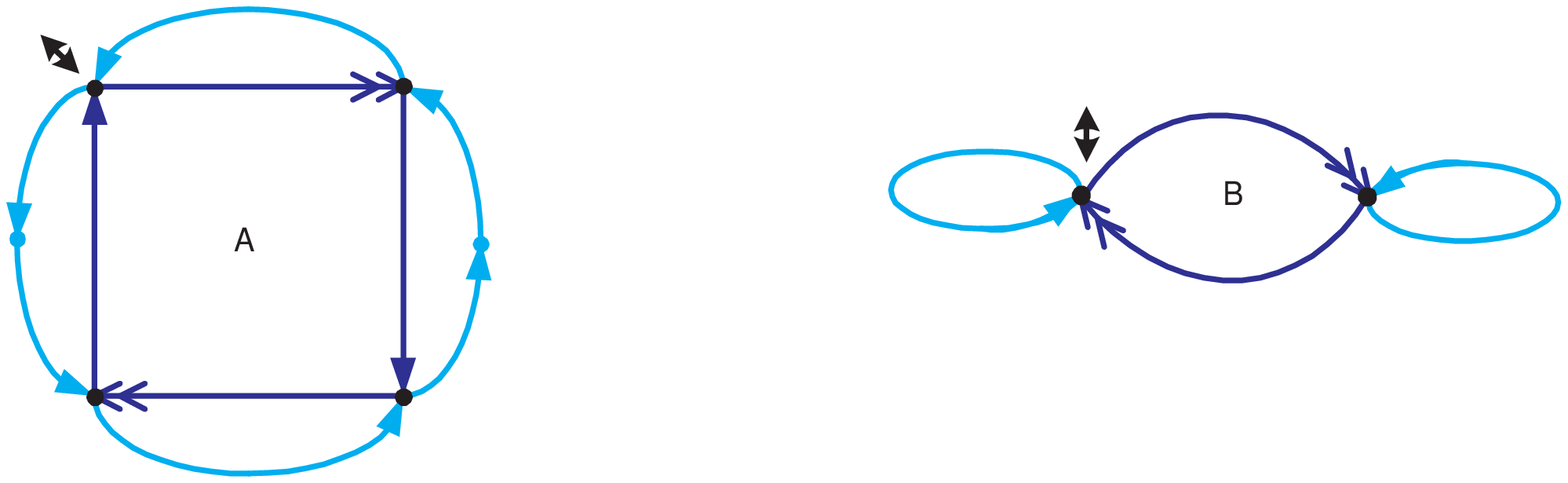}
\caption { \label{fig: ExOfFI}}
\end{figure}

\subsection*{Complexity.}
Let $m$ be the sum of the lengths of the words $h_1, \ldots h_k$.
By Theorem~\ref{thm: properties of subgroup graphs} (4), the
complexity of the construction of  $\Gamma(H)$ is $O(m^2)$.

The detecting of monochromatic components in this graph takes $ \:
O(|E(\Gamma(H))|) \: $. Since, by our assumption, all the
essential information about $A$, $G_1$ and  $G_2$ is given and it
is not a part of the input, verifications concerning a particular
monochromatic component of $\Gamma(H)$ take $O(1)$. Therefore to
do such verifications for all monochromatic component of
$\Gamma(H)$ takes $O(|E(\Gamma(H))|)$.
 Since, by Theorem~\ref{thm: properties of subgroup graphs} (4), $|E(\Gamma(H))|$ is proportional to $m$,  the
complexity of the ``freeness'' detecting  presented along with the
proof of Corollary \ref{freeness in amalgams of finite grp} is
$O(m^2)$, that is it is quadratic in the size of the input.

If the subgroup $H$ is given by the  graph $\Gamma(H)$, then to
verify that $H$ is a free subgroup of $G$ takes
$O(|E(\Gamma(H))|)$. That is the ``freeness'' algorithm  is even
linear in the size of the input.


\section{The Finite Index Problem}
\label{sec:FiniteIndexProblem}

One of the first natural computational questions regarding
subgroups is to compute the index of the subgroup in the given
group, the \emph{finite index problem}.

As is well known  (\cite{b-m-m-w, kap-m}), this problem is easily
solvable via subgroup graphs in the case of free groups. Recall
that a subgroup $H$ of a free group $FG(X)$ has  finite index if
and only if its subgroup graph $\Gamma_H$   constructed by the
Stallings' folding  algorithm is \emph{full}, i.e.
\emph{complete}, i.e. \emph{$X^{\pm}$-saturated}. That is for each
vertex $v \in V(\Gamma_H)$   and for each  $x \in X^{\pm}$ there
exists an edge which starts at $v$ and which is labelled by $x$.

In \cite{schupp} similar results were obtained for finitely
generated subgroups of certain Coxeter groups and surface groups
of an extra-large type.

In general, the index $[G:H]$ equals to the \emph{sheet number} of
the covering space, corresponding to the subgroup $H $, of the
standard $2$-complex representing the group  $G$ (\cite{stil}).
Thus if $G$ is finitely presented, the index of $H$ in $G$ is
finite if and only if the 1-skeleton of the corresponding covering
space is finite. That is  if and only if the relative Cayley
graph, $Cayley(G,H)$, is finite.

By Theorem~\ref{thm: properties of subgroup graphs} (2) and
Corollary~\ref{cor:PrecoversSubgrOfCayleyGr}, a subgroup graph
$(\Gamma(H),v_0)$ is a subgraph of $(Cayley(G,H), H \cdot 1)$.  It
turns out that there exists a strong connection between the index
of $H$ in $G$ and how ``saturated'' the graph $\Gamma(H)$ is. We
describe this connection in Theorem~\ref{fi} and use it to solve
the \emph{finite index problem} in amalgams of finite groups
(Corollary \ref{algorithm finite index_finite_grp}).

%
%
%
%

The complexity analysis of the presented algorithm is given at the
end of the section.

\begin{thm} \label{fi}
Let $H$ be a finitely generated subgroup of an amalgam of finite
groups $G=G_1 \ast G_2$.

Then $[G:H] < \infty$ if and only if $\Gamma(H)$ is
$X^{\pm}$-saturated.
\end{thm}
\begin{proof} The ``if'' direction is clear. Indeed,  if $\Gamma(H)$ is
$X^{\pm}$-saturated then, by Lemma~\ref{lemma1.5}, $\Gamma(H)$ is
isomorphic to $Cayley(G,H,H \cdot 1)$. Since, by Theorem~\ref{thm:
properties of subgroup graphs}, the graph $\Gamma(H)$ is finite,
$Cayley(G,H,H \cdot 1)$ is a finite graph. Hence
$[G:H]=|V(Cayley(G,H))| < \infty$.

To prove the opposite direction we assume that $\Gamma(H)$ is not
$X^{\pm}$-saturated. Note that since $\Gamma(H)$ is a precover of
$Cayley(G,H)$, any of its monochromatic component are either
$X^{\pm}_1$-saturated or $X^{\pm}_2$-saturated. Thus every
bichromatic vertex of $\Gamma$ is $X^{\pm}$-saturated and each
monochromatic vertex is either $X^{\pm}_1$-saturated or
$X^{\pm}_2$-saturated.

Let $v$ be a $X_1$-monochromatic vertex of $\Gamma$. Then by
Lemma~\ref{lemma2.12}, there is a path $p$ in normal form such
that $\iota(p)=v_0, \ \tau(p)=v$, and $w \equiv lab(p)$ is a word
in normal form.
%
%
Let $(w_1, \ldots, w_n)$ be a normal decomposition of $w$. Then
there is $x \in X_2 \setminus A$ \footnote{We assume that A is a
proper subgroup of $G_1$ and of $G_2$,  otherwise the amalgam $G_1
\ast_A G_1$ is a finite group and all computations are trivial in
our context}, such that $(w_1, \ldots, w_n,x)$ represents a word
$w' \in G$ in normal form. Now if $w_1 \in G_1$ (more precisely,
$w_1 \in G_1 \setminus A$, since $w$ is a normal word) or if $w_1
\in G_2$ but $xw_1 \in G_2 \setminus A$ then $(w')^n$ is in normal
form for all $n \geq 1$.

Otherwise $xw_1 \in G_2 \cap A$. Then there exists $y \in X_1
\setminus A$, such that $(w_1, \ldots, w_n,x,y)$ represents a word
$w'' \in G$ in normal form  and $(w'')^n$ is in normal form for
all $n \geq 1$. But neither $w'$ nor $w''$, and hence neither
$(w')^n$ nor $(w'')^n$  (for all $n \geq 1$) label a  path closed
at $v_0$ in $\Gamma(H)$.
Thus $(w')^n \not\in H$ and $(w'')^n \not\in H$, for all $n \geq
1$.

The existence of such elements shows that $H$ has infinite index
in $G$. Indeed, for all $n_1 \neq n_2$ and $g \in \{w',w''\}$ we
have $H(g)^{n_1} \neq H(g)^{n_2}$, because otherwise
$(g)^{n_1-n_2} \in H$. Thus, without loss of generality we can
assume that $n_1 > n_2$, then $n_1-n_2 \geq 1$ and we get a
contradiction.

\end{proof}

\begin{cor} \label{algorithm finite index_finite_grp}
Let $h_1, \ldots h_n \in G$.
Then there exists an algorithm which computes the index of the
subgroup $H=\langle h_1, \ldots, h_n \rangle$ in $G$.
\end{cor}
\begin{proof} We first construct the graph $\Gamma(H)$, using the
generalized Stallings' folding  algorithm.

Then we verify if this graph is $(X_1 \cup X_2)^{\pm}$-saturated.
If no, the subgroup $H$ has  infinite index in $G$, by
Theorem~\ref{fi}. Otherwise,  the index of $H$ in $G$ is finite
and $[G:H]=|V(\Gamma(H))|$.

\end{proof}

%
%
\subsection*{Complexity} \
Let $m$ be the sum of the lengths of the words $h_1, \ldots h_n$.
By Theorem~\ref{thm: properties of subgroup graphs} $(4)$, the
generalized Stallings' algorithm  computes $(\Gamma(H),v_0)$ in
time $O(m^2)$.
By the proof of  Corollary \ref{algorithm finite
index_finite_grp}, the detecting of the index takes time
proportional to $|E(\Gamma(H)|$. (Indeed, for each vertex of
$\Gamma(H)$ we have to check if it is bichromatic, which takes
$\sum_{v \in V(\Gamma(H))} deg(v)=2 |E(\Gamma(H)|$.)
Since, by Theorem~\ref{thm: properties of subgroup graphs} (4),
$|E(\Gamma(H)|=O(m)$, the complexity of the algorithm given along
with the proof of Corollary \ref{algorithm finite
index_finite_grp} is $O(m^2)$.

If the subgroup $H$ is given by $(\Gamma(H),v_0)$  and not by a
finite set of subgroup generators, then the above algorithm is
even linear in the size of the graph.


\begin{ex} \label{ex:FiniteIndex}
{\rm Let $H_1$ and $H_2$ be the subgroups considered in Example
\ref{ex:Freeness}.

Analyzing the ``saturation'' of the graphs  $\Gamma(H_1)$ and
$\Gamma(H_2)$ illustrated on Figure~\ref{fig: ExOfFI}, we see that
$[G:H_1]=\infty$, while $[G:H_2]=2$.} \e
\end{ex}


\section{The Separability Problem}
\label{sec:Separability}

A group $G$ is  \emph{subgroup separable}, or \emph{LERF}, if
given a  finitely generated subgroup $H$ of $G$ and  $g \not\in H$
there exists a finite index subgroup $K \leq G$ with $H \leq K$
and $g \not\in K$.
We call $K$ a \emph{separating subgroup}. If one places a topology
on $G$ (called the \emph{profinite topology} \cite{hall2}), by
taking the collection of finite index subgroups as a neighborhood
basis of 1, then $G$ is LERF if and only if all its finitely
generated subgroups are closed.

 \emph{LERF} was introduced by M.Hall \cite{hall1}, who
proved that free groups are LERF. This property  is preserved by
free products \cite{burns, ro}, but it is not preserved by direct
products: $F_2 \times F_2$ is not LERF \cite{a_g}. Free products
of LERF groups with finite amalgamation are LERF \cite{a_g}. In
general, the property is not preserved under free products with
infinite cyclic amalgamation \cite{long-niblo, ribs}. However
amalgams of free groups over a cyclic subgroup are LERF
\cite{b-b-s}, and, by \cite{gi_sep}, free products of a free group
and a LERF group amalgamated over a cyclic subgroup maximal in the
free factor are LERF as well.

Subgroup separability of some classes of hyperbolic groups was
widely exploited in papers of Gitik \cite{gi_sep, gi_doub,
gi_rips}. Long and Reid \cite{long-reid1, long-reid2} studied this
property in 3-manifold topology and in hyperbolic Coxeter groups.
Results on subgroup separability for right-angle Coxeter groups
and for Coxeter groups of extra-large type can be found in
\cite{gi_coxeter} and in \cite{schupp}, respectively. These papers
include  detailed algorithms which construct separating subgroups
using graph-theoretic methods.

\emph{M.Hall property} is closely connected with subgroup
separability. A groups $G$ is \emph{M.Hall} if and only if each of
its finitely generated subgroups is a free factor in a subgroup of
finite index in $G$.  M.Hall property of virtually free groups was
deeply studied in works of Bogopolskii \cite{b1, b2}, where a
criterion to determine whether a virtually free group is M.Hall
was given.


An algorithmic aspect of the LERF property can be formulated as
the \emph{separability problem}. It asks to find an algorithm
which constructs a separating subgroup $K$ for a given finitely
generated subgroup $H$ and $g \not\in H$.

Let us emphasize that the knowledge that $G$ has a solvable
decision problem does not provide yet an effective procedure to
solve this problem.  Thus, on the one hand, since  amalgams of
finite groups are LERF, by the result of Allenby and Gregorac
\cite{a_g}, the separability problem in this class of groups might
be  solvable. On the other hand, we are  interested to find an
efficient solution.
Below we adopt some ideas of Gitik introduced in \cite{gi_sep} to
develop such an algorithm (given along with the proof of the Main
Theorem  (Theorem~\ref{thm:MainTheoremSeparability})). Our main
result in this issue is summarized in the following theorem.
%
\begin{thm}[The Main Theorem] \label{thm:MainTheoremSeparability}
Let $G=G_1 \ast_A G_2$ be an amalgam of finite groups.

 The \underline{separability
problem}  for $G$  is solvable if one of the following holds
\begin{enumerate}
  \item[(1)] $A$ is cyclic,
 \item[(2)] $A$ is malnormal in at least one of the factors $G_1$ or
 $G_2$,
 \item[(3)] $A \leq Z(G_i)$, for some $i \in \{1,2\}$.
 \item[*]
 In particular, the separability problem is solvable if
 at least one of the factors ($G_1$ or $G_2$) is Abelian.
\end{enumerate}
\end{thm}

%

Recall that  given a finitely generated subgroup $H$ of an amalgam
of finite groups $G=G_1 \ast_A G_2$ the generalized Stallings'
algorithm constructs the canonical subgroup graph $\Gamma(H)$
which is a (reduced) precover of $G$ (Theorem~\ref{thm: properties
of subgroup graphs} (2)). Thus in order to prove our Main Theorem
we first show that each finite precover $(\Gamma,v_0)$ of $G$,
when $G$ satisfies one of the conditions $(1)-(3)$, can be
embedded in a finite $X_i$-saturated precover $(\Gamma',v_0)$ of
$G$ ($i \in \{1,2\}$). Then we prove that such a precover  can be
embedded in a finite cover $(\Gamma'',v_0)$ of $G$. Finally, we
take $K=Lab(\Gamma'',v_0)$ to be the separating subgroup. This
completes the proof of the Main Theorem.

Example~\ref{example: separability}  demonstrates the computation
of the separating subgroup $K$ for a given subgroup $H \leq G$.

{ \ } \\


The \emph{amalgam} of labelled graphs $\Gamma_1$ and $\Gamma_2$
along $\Gamma_0$ denoted by $\Gamma_1 \ast_{\: \Gamma_0}
\Gamma_2$, is the pushout of the following diagram in the category
of labelled graphs: 
$$\begin{array}{ccl}
  \Gamma_0 & \rightarrow & \Gamma_1  \\
  \downarrow & \searrow & \downarrow \\
  \Gamma_2 & \rightarrow & \Gamma_1 \ast_{\; \Gamma_0} \Gamma_2, \\
\end{array}$$
where $i_1:\Gamma_0  \rightarrow  \Gamma_1$ and $i_2:\Gamma_0
\rightarrow  \Gamma_2$ are injective maps and none of the graphs
need be connected. The amalgam depends on the maps $i_1$ and
$i_2$, but we omit reference to them, whenever it does not cause
confusion. It can be easily seen that amalgamation consists of
taking the disjoint union of graphs and performing the
identification prescribed by $i_1$ and $i_2$ and subsequent
\emph{foldings} (an identification of the terminal vertices of a
pair of edges with the same origin and the same label) until a
labelled graph is obtained \cite{gi_sep, stal}.


\begin{lem} \label{construction of X1-saturated precover}
Let  $\Gamma$ be a finite precover of an amalgamated free product
of finite groups $G=G_1 *_A G_2$. Then $\Gamma$ can be embedded in
a $X_1^{\pm}$-saturated precover of $G$ with finitely many
vertices.
\end{lem}
\begin{proof}
Any vertex  of a graph well-labelled with $X_1^{\pm} \cup
X_2^{\pm}$  has one of the following types:
\begin{itemize}
    \item It is bichromatic.
    \item It is $X_1$-monochromatic.
    \item It is $X_2$-monochromatic.
\end{itemize}
Since $\Gamma$ is a precover of $G$, the above types take the form
(respectively):
\begin{itemize}
    \item It is $X_1^{\pm} \cup X_2^{\pm}$-saturated.
    \item It is $X_1$-monochromatic and $X_1^{\pm}$-saturated.
    \item It is $X_2$-monochromatic and $X_2^{\pm}$-saturated.
\end{itemize}
The  proof is by induction on the number of vertices of the third
type. If no such vertices exist, then $\Gamma$ is already
$X_1^{\pm}$-saturated. Assume that $\Gamma$ has $m$
$X_2$-monochromatic vertices, and let $v$ be one of them.

Let $C$ be a $X_2$-monochromatic component, such that $v \in
VM_2(C)$. Let $S= A_v$ be the \emph{stabilizer of} $v$ by the
action of $A$ on the vertices of $C$, that is $A_v=\{ x \in A \; |
\; v \cdot x=v \} \leq A$, and let $A(v)=\{ v \cdot x \; | \; x
\in A \} \subseteq V(C)$ be the \emph{$A$-orbit} of $v$.

Consider $Cayley(G_1,S,S \cdot 1)$. Thus $A_{S \cdot 1}=S=A_v$ and
the $A$-orbit $A(S \cdot 1)=\{ (S \cdot 1) \cdot x  \; | \; x \in
A \}=\{ S x  \; | \; x \in A \} \subseteq  V(Cayley(G_1,S))$ is
isomorphic to  $A(v)$. Hence, taking $\Gamma_v= \Gamma \ast_{\{v
\cdot x=S  x \; | \; x \in A \}}Cayley(G_1,S)$, we get a finite
compatible graph whose monochromatic components are covers of the
factors $G_1$ or $G_2$. Therefore, by Corollary \ref{corol2.13},
$\Gamma_v$ is a precover of $\Gamma$.

Since $A_{S \cdot 1}=S=A_v$, the only identifications in
$\Gamma_v$ are between vertices  of  $A(v)$ and $A(S \cdot 1)$.
Since these are sets of monochromatic vertices of different
colors, no foldings are possible in $\Gamma_v$. Hence the graphs
$\Gamma$ and $Cayley(G_1, S)$ embed in $\Gamma_v$. Thus the images
in $\Gamma_v$ of the vertices of $A(v)$  (equivalently, of $A(S
\cdot 1)$) are bichromatic vertices, while the chromacity of the
images of other vertices of $\Gamma$ and  $Cayley(G_1, S)$ remains
unchanged. Hence
$$|VM_2(\Gamma_v)|=|VM_2(\Gamma)|-|A(v)|<m.$$

  Therefore $\Gamma_v$ is a finite precover of
$G$ with $|VM_2(\Gamma_v)|<m$ such that $\Gamma$ embeds in
$\Gamma_v$. This completes the inductive step.

\end{proof}

\begin{remark} \label{remark: separability X2-sat precover}
{\rm By the symmetric arguments if the conditions of
Lemma~\ref{construction of X1-saturated precover} hold then
$\Gamma$ can be embedded in a $X_i^{\pm}$-saturated precover of
$G$ ($i \in \{1,2\}$) with finitely many vertices.} \e
\end{remark}

The proof of Lemma~\ref{construction of X1-saturated precover}
yields the following technical result, which we employ later to
produce $X_i$-saturated precovers ($i \in \{1,2\}$).

\begin{cor} \label{cor:SeparabilityTechnical}
Let $G=G_1 \ast_A G_2$ be an amalgam of finite groups.

Let $\Gamma_i$   be a finite precover of $G$  (not necessary
connected) and let $v_i \in VM_i(\Gamma_i)$ ($i \in \{1,2\}$).

If   $A_{v_1}=A_{v_2}$ then $A(v_1)\simeq A(v_2)$, and
$\Gamma=\Gamma_1 \ast_{\{v_1 \cdot a=v_2 \cdot a \; | \; a \in
A\}}\Gamma_2$ is a finite precover of $G$ such that the graphs
$\Gamma_1$ and $\Gamma_2$ embed into the graph $\Gamma$.
\end{cor}


Now we consider $\Gamma$ to be a finite
$X_{\beta}^{\pm}$-saturated precover of $G$ ($\beta \in \{1,2\}$),
where $G=G_1 *_A G_2$ is an amalgamated free product of finite
groups. In the consequent lemmas,  it is  showen that if  one of
the conditions from Theorem~\ref{thm:MainTheoremSeparability} is
satisfied then $\Gamma$ can be embedded into a finite cover of
$G$.

Since the graph $\Gamma$ is  $X_{\beta}^{\pm}$, any vertex of
$\Gamma$ is either bichromatic or $X_{\beta}$-monochromatic.
Moreover, the graph $\Gamma$ is compatible, as a precover of $G$.
Hence any $A$-orbit consists of the vertices of the same type.
Therefore the set of $X_{\beta}$-monochromatic vertices of
$\Gamma$ can be viewed as a disjoint union of distinct $A$-orbits.
This enables us to consider the following notation.

For each $v \in VM_{\beta}(\Gamma)$ we set $n_v$ to be the number
of vertices in the $A$-orbit of $v$, that is $n_v=|A(v)|$. Recall
that  $$A_v \leq A, \ A(v) \simeq {A} / {A_v}, \ {\rm thus} \
|A|=|A(v)||A_v|.$$ 

Let $n(\Gamma)=\{ n_v | v \in VM_{\beta}(\Gamma) \}$. For each $n
\in n(\Gamma)$, assume that $\Gamma$ has $m$ different $A$-orbits,
each containing $n$ $X_{\beta}$-monochromatic vertices. Let
$\{v_i| 1 \leq i \leq m \}$ be the set of representatives of these
orbits. Denote $S_i= A_{v_i}$. Then for all $1 \leq i \leq m$, $|
S_i |=\frac{|A|}{n}$.

Assume that $A$ has $r$ distinct subgroups $S_j$ ($1 \leq j \leq
r$) of order $\frac{|A|}{n}$ and assume that $\Gamma$ has $m_j$
 representatives of distinct orbits $v_j \in VM_{\beta}(\Gamma)$ with
$A_{v_j}=S_j$. Hence $\sum^{r}_{j=1}{m_j}=m$.

With the above notation, we formulate Lemmas~\ref{separability:
A<=Z(G_i)}, \ref{separability: A is malnormal} and
\ref{separability: A is cyclic}.


\begin{lem} \label{separability: A<=Z(G_i)}
If $A$ is a center subgroup of $G_{\alpha}$ (that is $A \leq
Z(G_{\alpha})$ \footnote{Recall that the \emph{center of} $G$ is
the subgroup  \ $ Z(G)=\{g \in G \; | \; gx=xg, \; \forall x \in
G\}.$  }),
then any finite $X_{\beta}^{\pm}$-saturated precover of $G$  can
be embedded in a cover of $G$ with finitely many vertices ($1 \leq
\beta \neq \alpha \leq 2$).
\end{lem}
\begin{proof}
The proof is by induction on  $|n(\Gamma)|$.

Since $A \leq Z(G_{\alpha})$, $S_j$ is normal in $G_{\alpha}$ for
all $1 \leq j \leq r$. Therefore for each vertex $u \in
V(Cayley(G_{\alpha}, S_j))$, we have $A_u=S_j$. Indeed,
$$A_u=Lab(Cayley(G_{\alpha}, S_j),u) \cap A=g^{-1}S_jg \cap A=S_j \cap
A=S_j,$$ where $g \in G_{\alpha}$, such that $(S_j \cdot 1) \cdot
g = u$.
Thus distinct $A$-orbits of vertices in $Cayley(G_{\alpha},S_j)$
are isomorphic to each other and have length $n$. Their number is
equal to
$$\frac{|V(Cayley(G_{\alpha},S_j))|}{n}
=\frac{|G_{\alpha}|}{|S_j|} :
\frac{|A|}{|S_j|}=\frac{|G_{\alpha}|}{|A|}=[G_{\alpha}:A].$$

Let $t=[G_{\alpha}:A]$. Let $\Gamma_{1}$ be the disjoint union of
$t$ isomorphic copies of $\Gamma$ and let $\Gamma_{2}$ be  the
disjoint union  of $m_j$ isomorphic copies of
$Cayley(G_{\alpha},S_j)$, for all $1 \leq j \leq r$. Then both
$\Gamma_{1}$ and $\Gamma_{2}$ have $tm_j$ distinct isomorphic
$A$-orbits of length $n$.

Let $\{ w_{ji} \; | \; 1 \leq i \leq tm_j, \; 1 \leq j \leq r \}$
and $\{ u_{ji} \; | \; 1 \leq i \leq tm_j, \; 1 \leq j \leq r \}$
be the sets of representatives of these orbits in $\Gamma_{1}$ and
in $\Gamma_2$, respectively. Thus $A_{w_{ji}}=S_j=A_{u_{ji}}$, for
all $1 \leq i \leq tm_j$ and $1 \leq j \leq r$.
Let $\Gamma'$ be the amalgam of $\Gamma_{1}$ and $\Gamma_{2}$ over
these sets of vertices,
$$ \Gamma'=\Gamma_{1} \ast_{\{w_{ji} \cdot a=u_{ji} \cdot a \; | \;   a \in A  \}} \Gamma_{2}.$$
By Corollary~\ref{cor:SeparabilityTechnical}, $\Gamma'$ is a
finite precover of $G$ such that the graphs $\Gamma_1$ and
$\Gamma_2$ embed in it. Therefore the graph $\Gamma$ embeds in
$\Gamma'$ as well. Moreover, by construction, the graph $\Gamma'$
is  $X_{\beta}^{\pm}$-saturated,  and $n(\Gamma')=n(\Gamma)
\setminus \{n\}$. Thus $\Gamma'$ satisfies the inductive
assumption. We are done.

\end{proof}


\begin{lem} \label{separability: A is malnormal}
If $A$ is a malnormal subgroup of $G_{\alpha}$ then any finite
$X_{\beta}^{\pm}$-saturated precover of $G$  can be embedded in a
cover of $G$ with finitely many vertices ($1 \leq \beta \neq
\alpha \leq 2$).
\end{lem}
\begin{proof}
The proof is by induction on  $|n(\Gamma)|$.

Since $A$ is malnormal in $G_{\alpha}$,   for each vertex $u \in
V(Cayley(G_{\alpha}, S_j))$ ($1 \leq j \leq r$) such that $u=(S_j
\cdot 1) \cdot g$, where $g \in G_{\alpha} \setminus A$, we have
$$A_u=Lab(Cayley(G_{\alpha}, S_j),u) \cap A=g^{-1}S_jg \cap A=\{1\}.$$
Therefore $A(u)\simeq A$ and $|A(u)|=|A|$.  Thus
$V(Cayley(G_{\alpha}, S_j))$ form one $A$-orbit isomorphic to
$A(v_j)$ of length $n$ with $A_{S_j \cdot 1}=S_j=A_{v_j}$, and
$c=(|V(C)|-n)/|A|$ \ $A$-orbits isomorphic to $A(u)\simeq A$ of
length $|A|$ with, roughly speaking, a trivial $A$-stabilizer.

On the other hand, in $Cayley(G_{\beta})$ the number of distinct
$A$-orbits  of length $|A|$ with the trivial $A$-stabilizer is $d=
|V(Cayley(G_{\beta}) |/|A|=|G_{\beta} |/|A|=[G_{\beta}:A]$.

Let $\Gamma_1$ be the disjoint union of $d$ isomorphic copies of
$\Gamma$ and $c r$ isomorphic copies of $Cayley(G_{\beta})$. Let
$\Gamma_2$ be the union of disjoint unions of $m_j d$ isomorphic
copies of $Cayley(G_{\alpha},S_j)$, for all $1 \leq j \leq r$.
Then both $\Gamma_1$ and $\Gamma_2$ have $m_j d$ distinct
$A$-orbits of length $n$ isomorphic to $A(v_j)$, and $c d  r$
different isomorphic $A$-orbits of length $|A|$.

Let $\{ w_{ji} \; | \; 1 \leq i \leq m_j d, \; 1 \leq j \leq r \}$
and $\{ u_{ji} \; | \; 1 \leq i \leq m_j d, \; 1 \leq j \leq r \}$
be the sets of representatives of the orbits of length $n$ in
$\Gamma_{1}$ and in $\Gamma_2$, respectively. Hence
$A_{w_{ji}}=S_j=A_{u_{ji}}$, for all $1 \leq i \leq m_j  d$ and $1
\leq j \leq r$.

Let $\{ x_{l} \; | \; 1 \leq l \leq c d r \}$ and $\{ y_{l} \; |
\; 1 \leq l \leq c d r \}$ be the sets of representatives of the
orbits of length $|A|$ in $\Gamma_{1}$ and in $\Gamma_2$,
respectively. Then $A_{x_l}=\{1\}=A_{y_l}$, for all $1 \leq l \leq
c d r$.

Let $\Gamma'$ be the amalgam of $\Gamma_{1}$ and $\Gamma_{2}$ over
these sets of vertices,
$$ \Gamma'=\Gamma_{1} \ast_{\{w_{ji} \cdot a=u_{ji} \cdot a \; | \;   a \in A  \}
\cup \{x_l \cdot a=y_l \cdot a \; | \; a \in A \} } \Gamma_{2}.$$
By Corollary~\ref{cor:SeparabilityTechnical}, $\Gamma'$ is a
finite  precover of $G$ such that the graphs $\Gamma_1$ and
$\Gamma_2$ embed in it. Therefore the graph $\Gamma$ embeds in
$\Gamma'$ as well. Moreover, by construction, the graph $\Gamma'$
is  $X_{\beta}^{\pm}$-saturated, and $n(\Gamma')=n(\Gamma)
\setminus \{n\}$. Thus $\Gamma'$ satisfies the inductive
assumption. We are done.

\end{proof}


\begin{lem} \label{separability: A is cyclic}
If $A$ is cyclic then any finite $X_{\beta}^{\pm}$-saturated
precover of $G$ can be embedded in a cover of $G$ with finitely
many vertices ($\beta \in \{1,2\}$).
\end{lem}
\begin{proof}
Since $A$ is cyclic, $S_i=S_j$ for all $1 \leq i,j \leq m$, that
is $A_{v_i}=A_{v_j}$. Assume that $A_v=S$, for all $v \in \{v_i \:
| \: 1 \leq i \leq m \}$.

Consider $Cayley(G_{\alpha}, S)$ ($1 \leq \beta \neq \alpha \leq
2$).  For each vertex $u \in V(Cayley(G_{\alpha}, S))$, we have
$$A_u=Lab(Cayley(G_{\alpha}, S),u) \cap A=g^{-1}Sg \cap A,$$ where $g \in G_{\alpha}$, such that $(S \cdot 1) \cdot
g = u$. Thus $|A_u|  \leq |S|$ and therefore, since $A$ is cyclic,
$A_u \leq S \leq A$.


\begin{claim}
Let $ \alpha \in \{1,2\}$. 

Then there is $0 < N \in {\bf Z}$ such that $Cayley(G_{\alpha},
S)$ can be embedded into a finite $X_{\alpha}^{\pm}$-saturated
precover $C$ of $G$, whose $X_{\alpha} $-monochromatic vertices
form $N$ distinct $A$-orbits of length $n$  isomorphic to each
other, with the $A$-stabilizer $S$. More precisely,
$$VM_{\alpha}(C)=\bigcup_{i=1}^N A(v_i), \ {\rm such \ that } \
A_{v_i}=S \ ( \forall \ 1 \leq i \leq N).$$
\end{claim}
\begin{proof}[Proof of the Claim]
The proof is by induction on the number of prime factors of $|S|$.


Assume first that $|S|=p$   is prime.  By the above observation,
for all $u \in V(Cayley(G_{\alpha}, S))$, either $A_u=S$, $|A(u)|
=n$, or $A_u=\{1\}$, $|A(u)|=|A|$ (that is $A(u)\simeq A$).

Assume that $V(Cayley(G_{\alpha}, S))$ form $b$ distinct
isomorphic orbits of length $n$. Hence the number of distinct
$A$-orbits of $V(Cayley(G_{\alpha}, S))$ of length $|A|$
isomorphic to $A(u)$ with the trivial $A$-stabilizer is
$c=(|V(C)|-n \cdot b)/|A|$.

On the other hand, in $Cayley(G_{\beta})$ the number of distinct
$A$-orbits  of length $|A|$  with the trivial $A$-stabilizer is
$$d=\frac{|V(Cayley(G_{\beta}))|}{|A|}=\frac{|G_{\beta}|}{|A|}=[G_{\beta}:A].$$

Let $C_1$ be the disjoint union of $d$ isomorphic copies of
$Cayley(G_{\alpha}, S)$.  Let $C_2$ be the disjoint union of $c$
isomorphic copies of $Cayley(G_{\beta})$. Then both $C_1$ and
$C_2$ have $c d$ distinct $A$-orbits of length $|A|$.

Let $\{ x_{l} \; | \; 1 \leq l \leq c d \}$ and $\{ y_{l} \; | \;
1 \leq l \leq c d \}$ be the sets of representatives of these
orbits  in $C_{1}$ and in $C_2$, respectively. Then
$A_{x_l}=\{1\}=A_{y_l}$, for all $1 \leq l \leq c d $.

Let $C$ be the amalgam of $C_{1}$ and $C_{2}$ over these sets of
vertices,
$$ C=C_{1} \ast_{ \{x_l \cdot a=y_l \cdot a \; | \; a \in A \} } C_{2}.$$
By Corollary~\ref{cor:SeparabilityTechnical}, $C$ is a finite
precover of $G$ such that the graphs $C_1$ and $C_2$ embed in it.
Therefore the graph $Cayley(G_{\alpha}, S)$ embeds in $C$ as well.
Moreover, by  construction, the graph $C$ is
$X_{\alpha}^{\pm}$-saturated, and $VM_{\alpha}(C)$ form $N=b d$
distinct $A$-orbits  of length $n$ isomorphic to each other, with,
roughly speaking, an $A$-stabilizer $S$.


Assume now that $|S|$ is not a prime number. Let
$V(Cayley(G_{\alpha}, S))$ form $t_i$ distinct $A$-orbits of
length $\frac{|A|}{i}$ isomorphic to $A(u_i)$ with the
$A$-stabilizer $A_{u_i} \leq S$. Thus $|A_{u_i}|=i$, where $i \in
I=\{ i \: | \: 1 \leq i < |S|,  \ i \mid |S| \}$.

By the inductive assumption, $Cayley(G_{\beta},A_{u_i})$ can be
embedded into a finite $X_{\beta}^{\pm}$-saturated precover $C_i$
of $G$ whose $X_{\beta}$-monochromatic vertices form $k_i$
distinct $A$-orbits isomorphic to $A(u_i)$ of length
$\frac{|A|}{i}$ with the $A$-stabilizer $A_{u_i}$.

Let $l=lcm(\{k_i \; | \; i \in I \}$).
We take $C'_1$ be the disjoint union of $l$ isomorphic copies of
$Cayley(G_{\alpha}, S)$.  Let $C'_2$ be the union of disjoint
unions of $\left( t_i   \frac{l}{k_i} \right)$ isomorphic copies
of $C_i$, for all $i \in I$. Then both $C'_1$ and $C'_2$ have $t_i
l$ distinct  $A$-orbits of length $|i|$ isomorphic to $A(u_i)$.

Let $\{ w_{ij} \; | \; i \in I, \; 1 \leq j \leq t_i l \}$ and $\{
u_{ij} \; | \; i \in I, \; 1 \leq j \leq t_i l \}$ be the sets of
representatives of these orbits in $C'_{1}$ and in $C'_2$,
respectively. Hence $A_{w_{ij}}=A_{u_i}=A_{u_{ij}}$, for all $i
\in I$ and $1 \leq j \leq t_i l$.

Let $C$ be the amalgam of $C'_{1}$ and $C'_{2}$ over these sets of
vertices,
$$ C=C'_{1} \ast_{\{w_{ij} \cdot a=u_{ij} \cdot a \; | \;   a \in A  \}} C'_{2}.$$
By Corollary~\ref{cor:SeparabilityTechnical}, $C$ is a finite
precover of $G$ such that the graphs $C'_1$ and $C'_2$ embed in
it. Therefore the graph $Cayley(G_{\alpha}, S)$ embeds in $C$ as
well. Moreover, by  construction, the graph $C$ is
$X_{\alpha}^{\pm}$-saturated, and $VM_{\alpha}(C)$ form $N=t_n
  l$ distinct $A$-orbits  of length $n$ isomorphic to $A(v)$
with the $A$-stabilizer $S$. We are done.

\end{proof}

Let $\Gamma_1$ be the disjoint union of $N$ isomorphic copies of
$\Gamma$ and let $\Gamma_2$ be the  disjoint unions of $m$
isomorphic copies of $C$. Then both $\Gamma_1$ and $\Gamma_2$ have
$m N$ distinct  $A$-orbits of length $n$ isomorphic to $A(v)$ with
the $A$-stabilizer $S$. The standard arguments used in the proofs
of Lemmas~\ref{separability: A<=Z(G_i)} and \ref{separability: A
is malnormal} complete  the proof.

\end{proof}


\begin{proof}[Proof of the Main Theorem]

We first construct the graph $(\Gamma(H),v_0)$, using the
generalized Stallings' folding  algorithm.

Without loss of generality we can assume that $g$ is a normal
word. Since $g \not\in H$, then, by Theorem~\ref{thm: properties
of subgroup graphs} (3), $v_0 \cdot g \neq v_0$. Thus either $g$
is readable in $\Gamma(H)$, that is $v_0 \cdot g =v \in
V(\Gamma(H))$, or it is not readable.

Assume first that $v_0 \cdot g =v \in V(\Gamma(H))$.
We apply the algorithm described along with the proof of
Lemma~\ref{construction of X1-saturated precover} to embed  the
finite precover $Lab(\Gamma(H),v_0 )$  into a finite
$X_i^{\pm}$-saturated precover $(\Gamma,\vartheta)$, where
$\vartheta$ is the image of $v_0$, and we take $1 \leq i \neq j
\leq 2$, if $A$ is malnormal or central in $G_j$.

Now we embed $(\Gamma,\vartheta)$ into a finite cover of $G$,
using the appropriate algorithm given along with the proof of one
of Lemmas~\ref{separability: A<=Z(G_i)}, \ref{separability: A is
malnormal} or \ref{separability: A is cyclic}. Let $(\Phi, \nu)$
 be the resulting graph, where $\nu$ is the image of $\vartheta$.

Let $K=Lab(\Phi,\nu)$. By Theorem~\ref{fi}, $[G:K] < \infty$ and
$(\Phi,\nu)=(\Gamma(K),u_0)$. Since $$\Gamma(H) \subseteq \Gamma
\subseteq \Phi,$$ we have
$$Lab(\Gamma(H),v_0) \leq Lab(\Gamma,\vartheta) \leq Lab(\Phi,\nu).$$
Thus $H \leq K$. However $g \not\in K$, because the above graphs
are inclusions are embeddings. Therefore we are done.

Assume now that $g$ is not readable in $\Gamma(H)$. Let $g_1$ be
the longest prefix of $g$ that is readable in $\Gamma(H)$, that is
$v_0 \cdot g_1 =v \in V(\Gamma(H))$. Thus $v \in VM(\Gamma(H))$.
Without loss of generality, we can assume that $v \in
VM_1(\Gamma(H))$.

We glue to $\Gamma(H)$  a ``stem''  labelled by $g_2$ at $v$,
where $g \equiv g_1g_2$. Let $\Gamma$ be the resulting graph (see
Figure~\ref{fig:SeparabilityStem}).


\begin{figure}[!h]
\begin{center}
\psfrag{v0 }[][]{$v_0$}
\psfrag{v1 }[][]{$v_1$}
\psfrag{v }[][]{$v$}
\psfrag{g1 }[][]{$g_1$}
\psfrag{g2 }[][]{$g_2$}
\psfrag{g }[][]{$\gamma$}
\psfrag{u1 }[][]{$u_1$}
\psfrag{u2 }[][]{$u_2$}
\psfrag{um }[][]{$u_m$}
\psfrag{A }[][]{$\Gamma(H)$}
\psfrag{C }[][]{$C$}
\psfrag{D }[][]{$D$}
\includegraphics[width=0.8\textwidth]{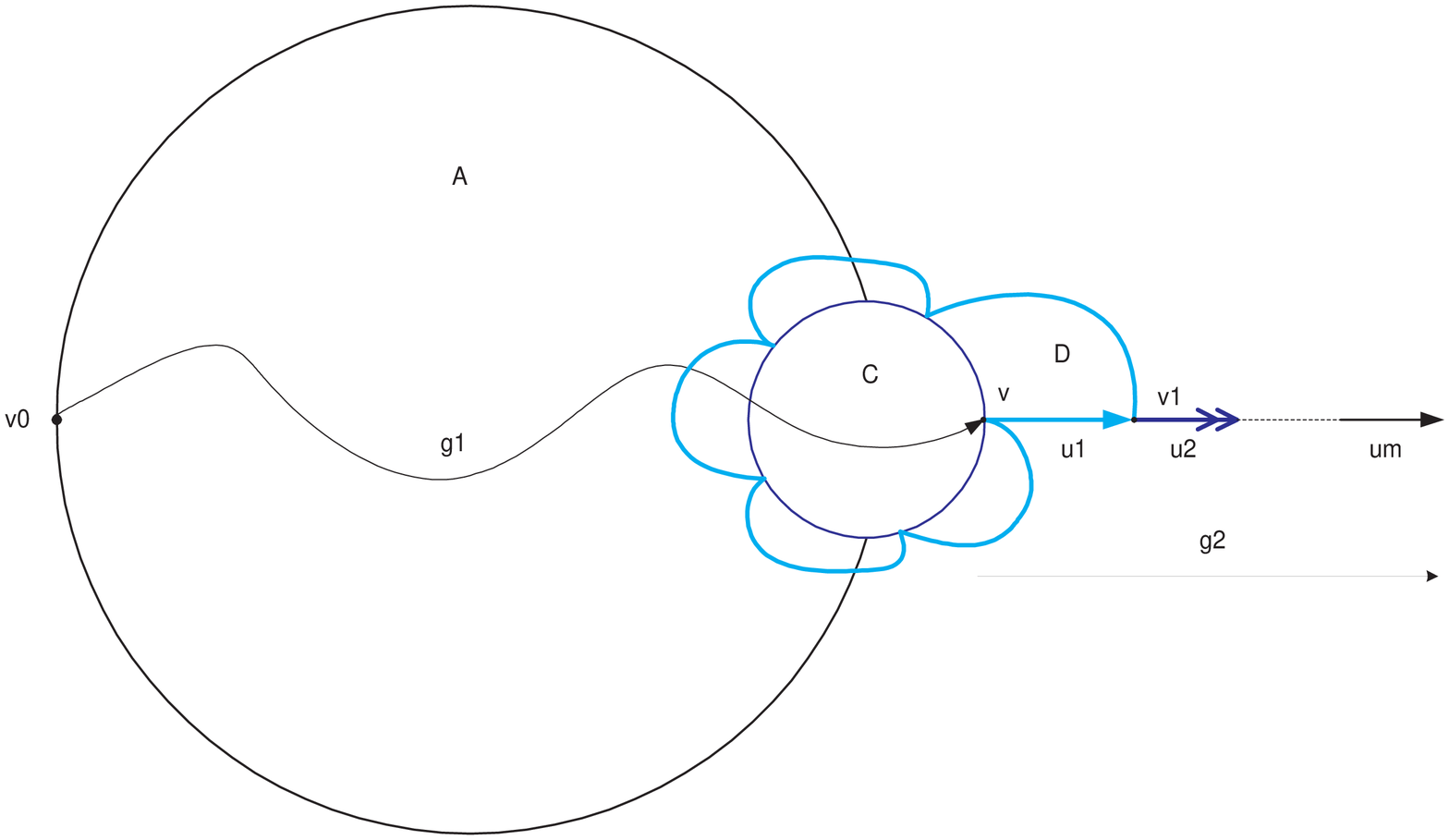}
\caption {  \label{fig:SeparabilityStem}}
\end{center}
\end{figure}


\begin{claim} \label{claim:PrecoverStem->Precover}
The graph $(\Gamma,v_0)$ can be embedded into a finite precover
$(\Gamma',v_0')$ of $G$ such that $v_0' \neq v_0' \cdot g \in
VM(\Gamma')$, where $v_0'$ is the image of $v_0$ in $\Gamma'$.
\end{claim}
\begin{proof}[Proof of the Claim]

Let $(u_1, \cdots, u_m)$ be the normal (Serre) decomposition of
$g_2$. Hence $u_1 \in G_2 \setminus A$. The proof is by induction
on the syllable length of $g_2$.

Let $C$ be a $X_1$-monochromatic component of $\Gamma(H)$ such
that $v \in V(C)$. Let $S=A_v$.

Consider $Cayley(G_1,S,S \cdot 1)$. Thus $A_{S \cdot 1}=S$ and the
$A$-orbit $A(S \cdot 1)=\{ (S \cdot 1) \cdot x  \; | \; x \in A
\}=\{ S x  \; | \; x \in A \} \subseteq  V(Cayley(G_2,S))$ is
isomorphic to the $A$-orbit of $v$ in $C$.
Therefore taking $\Gamma_v= \Gamma  \ast_{\{v \cdot x=S  x \; | \;
x \in A \}}Cayley(G_2,S)$, we get a graph such that $\Gamma(H)$
and $Cayley(G_2, S)$ embed in it, by
Corollary~\ref{cor:SeparabilityTechnical}.

Let $D$ be the $X_2$-monochromatic component of $\Gamma_v$ such
that $v \in V(D)$. Since $u_1 \in G_2$ and $D$ is
$X_2^{\pm}$-saturated, there exists a path $\gamma$ in $D$ such
that $\iota(\gamma)=v$ and $lab(\gamma) \equiv u_1$. Moreover, the
vertex $v_1=\tau(\gamma)  \in VB(D) \setminus VB(C)$, because $u_1
\in G_2 \setminus A$. Thus $v_1 \neq v_0$.

Therefore the graph $\Gamma_v$ can be thought of  as a precover of
$G$ with a stem labelled by $u_2 \cdots u_m$ which rises up from
the  vertex $v_1$. Note that $u_2 \in G_1 \setminus A$. Thus the
graph $\Gamma_v$ and the word given by the normal (Serre)
decomposition $(u_2, \cdots, u_m)$ satisfy the inductive
assumption. We are done.

\end{proof}


Proceeding in the same manner as in the previous case, when $v_0
\cdot g \in V(\Gamma(H))$, we embed the finite precover
$Lab(\Gamma',v_0' )$ of $G$ into a finite cover $(\Phi, \nu)$ of
$G$. This completes the proof.
 \end{proof}

\begin{figure}[!h]
\begin{center}
\psfrag{v0 }[][]{$v_0$}
\psfrag{v1 }[][]{$v_1$}
\psfrag{v2 }[][]{$v_2$}
\psfrag{v3 }[][]{$v_3$}
\psfrag{v4 }[][]{$v_4$}
\psfrag{v5 }[][]{$v_5$}
\includegraphics[width=0.8\textwidth]{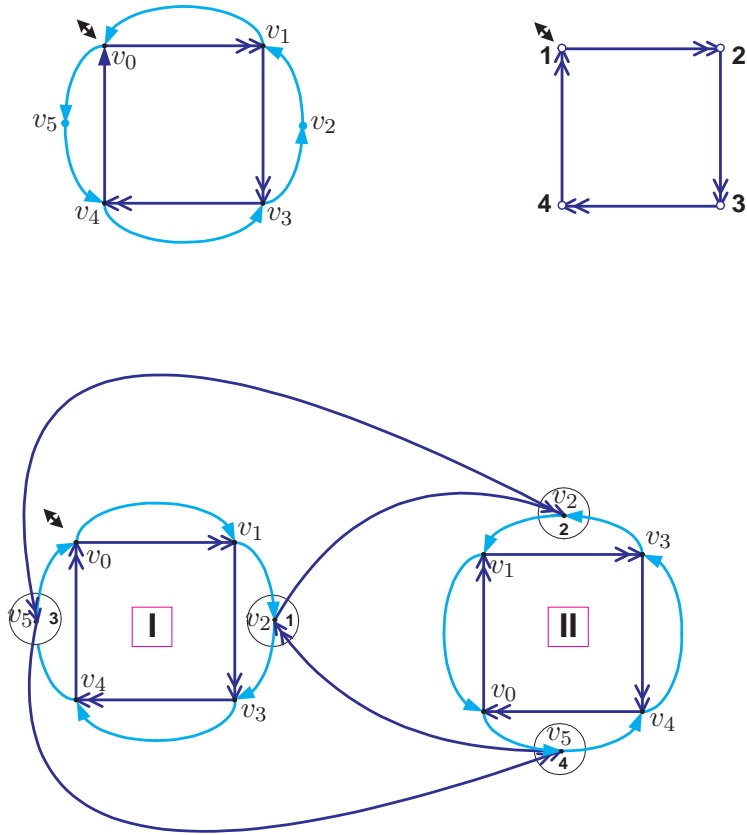}
\caption[The construction of the separating
subgroup]{{\footnotesize The construction of the cover $\Gamma(K)$
of $G$.} \label{fig: separability}}
\end{center}
\end{figure}

\begin{ex} \label{example: separability}
{\rm Let $G$ and $H_1$ be as in Example \ref{ex:Freeness}. Recall
that
$$G=\langle x,y | x^4, y^6, x^2=y^3 \rangle \ {\rm and} \ H_1=\langle xy \rangle.$$

Let $g=xy^{-1}$ be an element of $G$.  By Theorem~\ref{thm:
properties of subgroup graphs} (3), $g \not\in H_1$, because $v_0
\cdot g \neq v_0$.

Figure \ref{fig: separability} illustrates the construction of the
cover  $\Gamma(K)$ of $G$, where $K \leq G$ is the separating
subgroup for $H_1 \leq G$ and the element $g \not\in H_1$.}
 \e

\end{ex}


\appendix
\section{}


Below we follow the notation of Grunschlag \cite{grunschlag},
distinguishing between the ``\emph{input}'' and the ``\emph{given
data}'', the information that can be used by the algorithm
\emph{``for free''}, that is it does not affect the complexity
issues.

\begin{center}
\large{\emph{\underline{\textbf{Algorithm}}}}
\end{center}

\begin{description}
\item[Given] Finite groups $G_1$, $G_2$, $A$ and the amalgam
$G=G_1 \ast_{A} G_2$ given via $(1.a)$, $(1.b)$ and $(1.c)$,
respectively.

We assume that the Cayley graphs and all the relative Cayley
graphs of the free factors are given.
\item[Input]  A finite set $\{ g_1, \cdots, g_n \} \subseteq G$.
\item[Output] A finite graph $\Gamma(H)$ with a basepoint $v_0$
which is a reduced precover of $G$ and the following holds
\begin{itemize}
 \item
$Lab(\Gamma(H),v_0)=_{G} H$;
 \item $H=\langle g_1, \cdots, g_n \rangle$;
 \item a normal word $w$ is in $H$ if and only if
  there is a loop (at $v_0$) in $\Gamma(H)$
labelled by the word $w$.
 \end{itemize}
\item[Notation] $\Gamma_i$ is the graph obtained after the
execution of the $i$-th step.

%
%
\medskip

    \item[\underline{Step1}] Construct a based set of $n$ loops around a common distinguished
vertex $v_0$, each labelled by a generator of $H$;
    \item[\underline{Step2}] Iteratively fold edges and cut hairs %
    \footnote{A \emph{hair} is an edge one of whose endpoint has degree 1};
  \item[\underline{Step3}] { \ }\\
\texttt{For} { \ } each $X_i$-monochromatic component $C$ of
$\Gamma_2$ ($i=1,2$) { \ } \texttt{Do} \\
\texttt{Begin}\\
    pick an edge $e \in E(C)$; \\
    glue a copy  of $Cayley(G_i)$   on $e$ via identifying $ 1_{G_i} $  with $\iota(e)$ \\
    and identifying the two copies of $e$ in $Cayley(G_i)$ and in $\Gamma_2$; \\
    \texttt{If}  { \ } necessary  { \ } \texttt{Then} { \ } iteratively fold
    edges; \\
 \texttt{End;}

 \item[\underline{Step4}] { \ } \\
\texttt{ For}  { \ } each $v \in VB(\Gamma_3)$ { \ } \texttt{ Do} \\
 \texttt{If} { \ } there are paths $p_1$ and $p_2$, with $\iota(p_1)=\iota(p_2)=v$
 and $\tau(p_1)~\neq~\tau(p_2)$  such that
 $$lab(p_i) \in G_i \cap A \ (i=1,2) \ {\rm and} \  lab(p_1)=_G
 lab(p_2)$$
\texttt{ Then} { \ } identify $\tau(p_1)$ with $\tau(p_2)$; \\
 \texttt{If}  { \ } necessary  { \ } \texttt{Then} { \ } iteratively fold
    edges; \\

 \item[\underline{Step5}]
%
Reduce  $\Gamma_4$ by an iterative removal of all
(\emph{redundant})
 $X_i$-monochromatic components $C$ such that
\begin{itemize}
 \item $(C,\vartheta)$ is isomorphic to $Cayley(G_i, K, K \cdot 1)$, where $K \leq A$ and
$\vartheta \in VB(C)$;
 \item  $|VB(C)|=[A:K]$;
 \item one of the following holds
    \begin{itemize}
        \item  $K=\{1\}$ and $v_0 \not\in VM_i(C)$;
        \item $K$ is  a nontrivial subgroup of $A$ and $v_0  \not\in V(C)$.\\
    \end{itemize}
\end{itemize}

Let $\Gamma$ be the resulting graph;\\

\texttt{If}  { \ }
$VB(\Gamma)=\emptyset$ and $(\Gamma,v_0)$ is isomorphic to $Cayley(G_i, 1_{G_i})$ \\
\texttt{Then} { \ } we set $V(\Gamma_5)=\{v_0\}$ and
$E(\Gamma_5)=\emptyset$; \\
\texttt{Else} { \ } we set $\Gamma_5=\Gamma$.

 \item[\underline{Step6}] { \ } \\
 \texttt{If} { \ }
 \begin{itemize}
  \item $v_0 \in VM_i(\Gamma_5)$ ($i \in \{1,2\}$);
  \item $(C,v_0)$ is isomorphic to $Cayley(G_i,K,K \cdot 1)$, where $L=K \cap A$ is a nontrivial
  subgroup of
 $A$ and $C$ is a $X_i$-monochromatic component of $\Gamma_5$ such that $v_0 \in V(C)$;
  \end{itemize}
\texttt{Then} { \ } glue to $\Gamma_5$ a $X_j$-monochromatic
component ($1 \leq i \neq j \leq 2$) $D=Cayley(G_j,L,L \cdot 1)$
via identifying $L \cdot 1$ with $v_0$ and \\
identifying the vertices $L \cdot a$ of  $Cayley(G_j,L,L \cdot 1)$
with the vertices $v_0 \cdot a$ of $C$, for all $a \in A \setminus
L$.

Denote $\Gamma(H)=\Gamma_6$.

\end{description}


\begin{remark} \label{stal-mar-meak-kap-m}
{\rm Note that the first two steps of the above algorithm
correspond precisely to the Stallings' folding algorithm for
finitely generated subgroups of free groups \cite{stal, mar_meak,
kap-m}.} \e
\end{remark}


\begin{figure}[!h]
\psfrag{x }[][]{$x$} \psfrag{y }[][]{$y$} \psfrag{v }[][]{$v$}
\psfrag{x1 - monochromatic vertex }[][]{{\footnotesize
$\{x\}$-monochromatic vertex}}
\psfrag{y1 - monochromatic vertex }[][]{\footnotesize
{$\{y\}$-monochromatic vertex}}
\psfrag{ bichromatic vertex }[][]{\footnotesize {bichromatic
vertex}}
\includegraphics[width=0.8\textwidth]{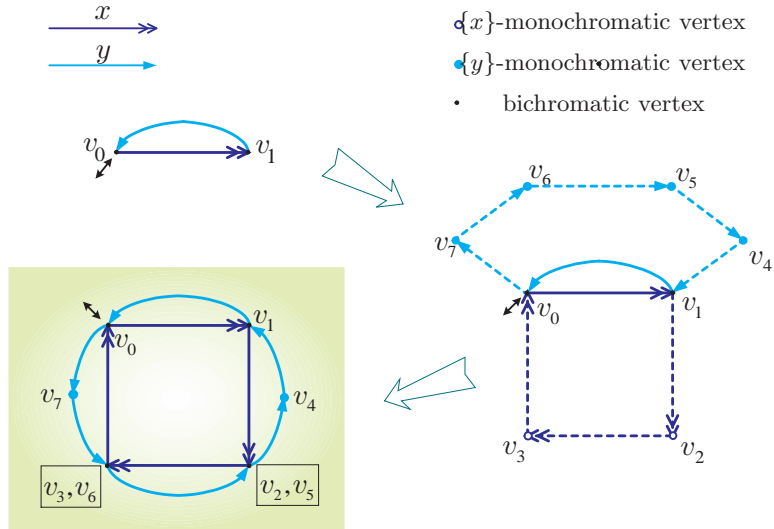}
\caption[The construction of $\Gamma(H_1)$]{ \footnotesize {The
construction of $\Gamma(H_1)$.}
 \label{fig: example of H=xy}}
\end{figure}

\begin{figure}[!hb]
\psfrag{v }[][]{$v$}
\includegraphics[width=0.8\textwidth]{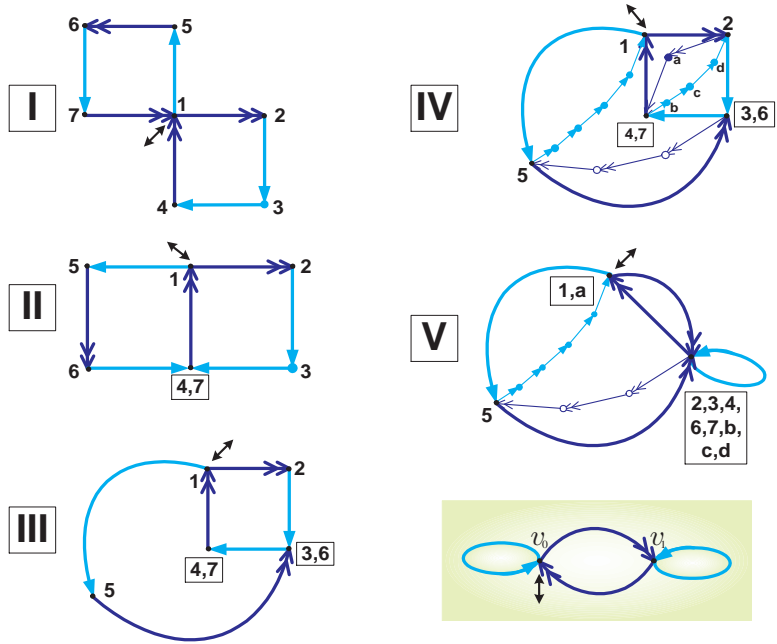}
\caption[The construction of $\Gamma(H_2)$]{ \footnotesize {The
construction of $\Gamma(H_2)$.} \label{fig: example of H=xy^2x,
yxyx}}
\end{figure}

\begin{ex} \label{example: graphconstruction}
{\rm Let $G=gp\langle x,y | x^4, y^6, x^2=y^3 \rangle$.

Let $H_1$ and $H_2$ be  finitely generated subgroups of $G$ such
that
$$H_1=\langle xy \rangle \ {\rm and} \ H_2=\langle xy^2, yxyx \rangle.$$

The construction of $\Gamma(H_1)$ and $\Gamma(H_2)$ by the
algorithm presented above is illustrated on Figures \ref{fig:
example of H=xy}
 and  \ref{fig: example of H=xy^2x, yxyx}.}
\e
\end{ex}


\end{document}